\documentclass[11pt]{article}
\usepackage{epsfig}
\usepackage{amsmath}
\usepackage{amsfonts}
\usepackage{amssymb}
\usepackage{color}

\title{Bifurcations of stationary measures of random diffeomorphisms}

\author{Hicham Zmarrou and Ale Jan Homburg\\
KdV Institute for Mathematics \\
University of Amsterdam \\
Plantage Muidergracht 24 \\
1018 TV Amsterdam \\
The Netherlands}

\newcounter{bean}

\newtheorem{theorem}{Theorem}[section]

\newtheorem{lemma}[theorem]{Lemma}
\newtheorem{proposition}[theorem]{Proposition}
\newtheorem{definition}[theorem]{Definition}
\newtheorem{remark}[theorem]{Remark}

\newcommand{\Res}[2]{{#1}\raisebox{-.4ex}{$\left|\,_{#2}\right.$}}
\newcommand{\qed}{\hspace*{\fill}$\square$}

\newcounter{hyp}
\newtheorem{hypothesis}{}

\newenvironment{hyp}{ \stepcounter{hyp}
    \begin{hypothesis} \rm \rule{2ex}{0pt}
    \begin{minipage}[t]{12.8cm} \rule{0pt}{0pt}\\[-1.25em]
\rule{0pt}{0pt}\!\!}{%
    \end{minipage}
   \end{hypothesis}}

\begin{document}
\maketitle

\footnote{2000 Mathematics Subject Classification  37A50,37Hxx,37Gxx,60Gxx}

\begin{abstract}
Random diffeomorphisms with bounded absolutely continuous noise
are known to possess a finite number of stationary measures.
We discuss dependence of stationary measures on an auxiliary parameter, thus
describing bifurcations of families of random diffeomorphisms.
A bifurcation theory is developed under mild regularity assumptions 
on the diffeomorphisms and the noise
distribution (e.g. smooth diffeomorphisms with uniformly distributed additive
noise are included).
We distinguish bifurcations where the density function 
of a stationary measure varies discontinuously or
where the support of a stationary measure varies discontinuously.

We establish that generic random diffeomorphisms are stable.
Densities of stable stationary measures are shown
to be smooth and
to depend smoothly on an auxiliary parameter,
except at bifurcation values.
The bifurcation theory explains the occurrence
of transients and intermittency as the main bifurcation phenomena in random diffeomorphisms.
Quantitative descriptions by means of average escape times
from sets as functions of the parameter are provided.
Further quantitative properties are described through
the speed of decay of correlations as function of the parameter.

Random endomorphisms are studied in one dimension; we show
that stable one dimensional random endomorphisms occur open and dense and that
in  one parameter families bifurcations are typically isolated.
We classify codimension one bifurcations for one
dimensional random endomorphisms;
we distinguish
three possible kinds, 
the random saddle node, the random homoclinic and the
random boundary bifurcation.
The theory is illustrated on families of random circle diffeomorphisms
and random unimodal maps.
\end{abstract}

\newpage
\tableofcontents

\section{Introduction}

To fix thoughts, consider a single map with an attracting fixed point $P$.
Write $W^s(P)$ for its basin of attraction.
Adding uniform noise of small amplitude gives a random map with
a stationary density with support near $P$.
Increasing the amplitude of the noise leads to a bifurcation when orbits
can escape from $W^s(P)$.
Two possibilities occur: escaping orbits can or cannot return near $P$.
Escaping orbits lead to transient dynamics if orbits do not return or intermittent dynamics
if orbits do return.
How can such transitions occur?
What are the quantitative characteristics?
As a second issue, take a one parameter family of maps that
exhibits bifurcations and add small bounded noise to it.
What happens to the bifurcation set and
in what way do bifurcations in the randomly perturbed map
manifest?
More generally, one can consider maps and families of maps where noise is an
intrinsic part of the description.

It is the purpose of this paper to work out bifurcation theory
of random smooth diffeomorphisms from the perspective of
stationary measures, providing answers to questions like the
ones just stated.
The context in which we perform this study is that of
points being mapped into bounded domains according to a
probability distribution, under some regularity conditions.
One can think of points being mapped by a diffeomorphism, defined 
on some compact manifold $\mathcal{M}$, followed by a
random perturbation.
This defines a discrete Markov process on  $\mathcal{M}$
given by transition functions $P(x,A)$
providing the chance that a point $x\in \mathcal{M}$ ends up in a Borel set $A \subset \mathcal{M}$.
We will assume that the region where $x$ is mapped into, is a bounded
domain $U_x$.
We argue that from a modeling point of view there are clear and
good reasons to consider bounded noise:
in most physical systems random perturbations are limited in their effect.

We indicate the regularity conditions assumed in this paper. 
A first regularity assumption is smoothness of the density
of the transition functions $P(x,A)$. It is not assumed that
this density vanishes on the boundary of its support; densities that are positive on bounded domains to model uniform noise are incorporated.
For $y \in \mathcal{M}$,
write $V_y$ for the set of points $x \in \mathcal{M}$ that are mapped to domains that include $y$;
$$
V_y = \{ x \in \mathcal{M} \;|\; y \in \text{support } P(x,\cdot)\}.
$$
As a second regularity assumption, we suppose that the sets $V_y$ are
domains with piecewise smooth boundary varying
smooth with $y$.
This assumption is natural in the context of
diffeomorphisms followed by noise, but does not hold in the context
of endomorphisms (possessing critical points) followed by noise.
Precise formulations follow in Section~\ref{subsec_maps} below.

An alternative description of the setup is by starting with
a collection ${\cal F}$ of maps
on $\mathcal{M}$ and a measure on ${\cal F}$. Maps are drawn randomly, and independently,
from ${\cal F}$ according to the given measure.
Similarly one can consider
maps depending on parameters that are drawn randomly.
The random parameters are drawn from a bounded domain according to a given distribution.
Typical examples are given by smooth maps $f$
with additive or parametric noise.
It will turn out that there is no loss of generality, as far as
statistical descriptions are concerned, when considering
maps with finitely many random parameters.
We address an appendix to the exploration of the range and connections of these
definitions.

In the following we will mostly speak of random maps;
maps depending on finitely many parameters that are drawn
randomly.
In the part of the paper developing the general theory,
Sections~\ref{sec_transfer} to \ref{sec_correlations}, we
assume the maps to be diffeomorphisms.
This guarantees
the regularity assumptions formulated above in the description
as discrete Markov processes to hold.
Random endomorphisms are studied in one dimension in Section~\ref{sec_1D}.

Under the mild regularity assumptions,
random diffeomorphisms possess finitely many
stationary measures, whose support is the closure of an open set and
whose density functions (stationary densities)
are smooth on all of $\mathcal{M}$. The
stationary densities are flat along the boundary of their support.
This has immediate consequences for the statistical properties of orbits:
orbits are very rarely found near the boundary of the
support of the stationary densities.
The main focus of this paper lies then
in the description of the dependence
of stationary densities on the random diffeomorphisms.
This includes describing quantitative characteristics.
Motivated by examples discussed below we call a random diffeomorphism stable
if its stationary densities and their supports vary continuously
with the random diffeomorphism (precise definitions follow shortly).
Otherwise we speak of a bifurcation.

Below we introduce the precise setup and present
our main results in a series of theorems.
The presentation of the material is
 separated in a section treating random diffeomorphisms
(Section~\ref{subsec_maps})
and a section treating families of random
diffeomorphisms (Section~\ref{subsec_families}).
Various aspects of the bifurcation theory for random diffeomorphisms
are presented in Theorems~\ref{stationary_measures},
\ref{isolated}, \ref{stable}, \ref{continuation}, \ref{theorem_exittimes_smooth} and \ref{theorem_corr}.
Appendix~\ref{sec_representation}
comments on the setup.
The main body of the theory is developed in Sections~\ref{sec_transfer} to \ref{sec_correlations}.
Sections~\ref{sec_transfer}, \ref{sec_stable} and \ref{sec_parameter}
contain respectively
material on transfer operators, proofs of stability theorems,
and discussions of parameter dependence.
Sections~\ref{sec_escape} and \ref{sec_escapeII} develop material on conditionally
stationary measures and apply this to compute expected escape times.
Section~\ref{sec_correlations} treats the speed of decay of correlations
depending on a parameter.
Appendix~\ref{sec_ift} contains an implicit function theorem used
to obtain regularity in the parameter of solutions of integral equations
involving the transfer operator.

There is no easy analogous theory for random endomorphisms;
stationary measures for random endomorphisms will in general
be less regular resulting in different statistical properties.
In Section~\ref{sec_1D} we present
a satisfactory theory for one dimensional random endomorphisms,
including a classification of possible bifurcations.
This extends the general theory
in particular by classifying codimension one bifurcations.
Typically only finitely many bifurcations occur in one parameter families
of random endomorphisms, in contrast to families of deterministic maps.
The material on one dimensional random endomorphisms
is developed and presented in Section~\ref{sec_1D}.

Section~\ref{sec_casestudies} contains two worked out examples
of random circle diffeomorphisms and random unimodal maps.
In Section~\ref{sec_circle} we consider
the standard circle diffeomorphism with
small additive noise;
\begin{equation}
f_a (x;\omega) = x + a + \omega + \varepsilon \sin(2 \pi x)
\end{equation}
with $x$ on the circle ${\mathbb R}/{\mathbb Z}$,
for fixed $\varepsilon \in (0,1)$ and a random parameter
$\omega$ taken uniformly from a small interval.
Section~\ref{sec_logistic} discusses as a prototypical example
of random endomorphisms on an interval,
random logistic maps
\begin{equation}
f (x;\mu) = (\mu + \omega)  x (1-x)
\end{equation}
with multiplicative noise obtained by varying $\omega$
with a uniform distribution in some interval.

There is a large body of literature
on stochastic stability (e.g. \cite{you86,kif88,balvia96,alvara03,alvaravas04},
see also \cite{bondiavia05}),
considering bounded noise as a means to treat properties of
single deterministic systems. This is done by letting the noise level decrease to zero.
In contrast, we consider maps and families of maps
where noise is an intrinsic part of the description.

Previous attempts to study stochastic bifurcations fall into two categories.
One way is to consider
notions close to the traditional understanding of bifurcations
for deterministic dynamics by
embedding the random dynamical system into a skew product system.
Another way studied frequently in the literature uses
singularity theory to describe changes in the density of a stationary measure.
This approach has been used to study systems with unbounded noise, often of a Gaussian nature,
so that a unique smooth stationary density occurs.
See \cite{arn98} for further discussion.
We consider the shape of the stationary densities, but the noise being bounded
leaves dynamics part of the picture.

Central in control theory are control sets, that is
maximal sets of approximate controllability.
The shape of the control sets varies with external
parameters, where discontinuous changes are possible.
Interpreting the control variable as noise, a relation with
the present paper becomes apparent.
Discontinuous changes of the control set, interpreted
this way, are among the
bifurcations identified in this paper.
In the context of differential equations depending on
a control variable the study of control sets is taken up in
\cite{colkli00,gay04}.
Bifurcation theory for such random differential equations 
in the spirit of this paper is considered in \cite{homyou05b}.

\subsection{Random diffeomorphisms}\label{subsec_maps}

The adjective smooth stands for $C^\infty$.
Let $\mathcal{M}$ be a smooth $n$-dimensional compact Riemannian manifold
with measure $m$ induced by the Riemannian structure.
Let $\Delta$ be a closed domain in $n$-dimensional Euclidean space.
Smoothness of a function $g$ on $\Delta$ is to be understood
in the sense that $g$ can be extended
to a smooth function on a neighborhood of $\Delta$.

\begin{definition}
A {\em smooth random map}, or {\em random endomorphism}, is a
smooth map $f: \mathcal{M} \times \Delta \to \mathcal{M}$, $x \mapsto f(x;\omega)$,
depending on a random parameter $\omega \in \Delta$
drawn from a measure on $\Delta$ with
smooth density function
$g: \Delta \to {\mathbb R}$, $\omega \mapsto g(\omega)$.
A {\em random diffeomorphism} is a smooth random map so that
$x \mapsto f(x;\omega)$ is a diffeomorphism for each $\omega$.
\end{definition}

\begin{remark}\label{remark_g}
Alternatively one can explicitly include the noise distribution $g$
as part of the definition of smooth random map (and speak of a pair
$(f,g)$). For convenience we consider $g$ given. The results in this paper
have direct, easily obtained, analogs if $g$ is allowed to vary.
\end{remark}

Note that endomorphisms and diffeomorphisms are always assumed to be smooth.
The basic setup we are treating is
of points
being mapped into bounded domains according to some
probability. The following standing assumptions will be made with this setup in mind.
The random parameters will be chosen from a region $\Delta$ that is a domain
in ${\mathbb R}^n$
with a piecewise smooth boundary.
The number of random parameters is in particular equal to the dimension of
the state space $\mathcal{M}$.
The most important examples are
where $\Delta$ is the unit ball
$\Delta = \{ x \in {\mathbb R}^n \;|\; \|x\| \le 1\}$
or the unit box
$\Delta = \{ (x_1,\ldots,x_n) \in {\mathbb R}^n \;|\; |x_1|,\ldots,|x_n| \le 1\}$.
Throughout this paper we assume that $\omega \mapsto f(x;\omega)$
is an injective map for each $x$.
Hence $f(x; \Delta)$ is diffeomorphic to $\Delta$.

Write $\nu$ for the measure on $\Delta$ with density function $g$.
A smooth random map gives rise to a discrete Markov process through
the transition functions
\begin{equation}\label{transition}
P(x,A) =  \int_{\{ \omega \;|\;  f(x;\omega) \in A\}} d\nu(\omega)
\end{equation}
for Borel sets $A$.
With $h_x(\omega) = f(x;\omega)$,
the measure $P(x,\cdot)$ equals $(h_x)_* \nu$ defined
by $(h_x)_* \nu (A) = \nu( h_x^{-1} (A))$.
Vice versa, a discrete Markov process with noise from a ball or a box
such that its transition functions
have smooth positive densities admits a representation by smooth random maps
(depending injectively on a random parameter),
see Appendix~\ref{sec_representation}.
Some cases of parametric noise, where the maps
do not depend
injectively on the random parameter, do not fall into this setup.
We refer to Appendix~\ref{sec_representation} for further discussion.

The general theory will be developed for random diffeomorphisms, instead of
random endomorphisms.
Endomorphisms allow for pathological
examples, for instance maps $f(x;\omega)$ that are constant in $x$.
We will however discuss random endomorphisms
in one dimension (on a circle or a compact interval) in detail.

With a slight abuse of notation,
iterates of $f(x;\omega)$ are given as
\begin{equation}\label{iterate}
f^k(x; \omega_1,\ldots,\omega_k) =
   f (f^{k-1}(x;\omega_1,\ldots,\omega_{k-1});\omega_k).
\end{equation}
More generally, write $\Delta^{\mathbb N}$ for all infinite
sequences ${\boldsymbol \omega} = \{\omega_i\}_{i \ge 1}$ with
each $\omega_i \in \Delta$. Denote $f^k(x;{\boldsymbol \omega}) =
f^k(x; \omega_1,\ldots,\omega_k)$. Let $\vartheta :
\Delta^{\mathbb N} \to \Delta^{\mathbb N}$ be the left shift
operator; $\vartheta \{\omega_i\}_{i \ge 1} = \{\omega_i\}_{i \ge
2}$. Consider the skew product system $S: \mathcal{M} \times \Delta^{\mathbb
N} \to \mathcal{M} \times \Delta^{\mathbb N}$ given by
\begin{equation}\label{skew}
S(x,{\boldsymbol \omega}) = (f(x;\omega_1) , \vartheta
{\boldsymbol \omega}).
\end{equation}
On $\Delta^{\mathbb N}$ one considers a measure $\nu^\infty$
which is the product of the measure $\nu$ over each $\Delta$.

With these definitions in mind, we introduce the central notions
of stationary measures and ergodic measures.
A stationary measure $\mu$ for the smooth random map $f$
is a probability measure on $\mathcal{M}$ with
$\mu \times \nu^\infty$ S-invariant;
$$
\mu \times \nu^\infty (S^{-1} (B)) = \mu \times \nu^\infty (B)
$$
for Borel sets $B \subset \mathcal{M} \times \Delta^{\mathbb N}$.
Equivalently, see \cite{kif88,ara00},
$$
\mu(A) = \int_\mathcal{M} P(x,A) d\mu (x)
$$
for Borel sets $A\subset \mathcal{M}$.
We refer to a stationary density as the density of an absolutely continuous
stationary measure.

A stationary measure $\mu$ is called ergodic if $\mu\times \nu^\infty$
is an ergodic measure for $S$ in the usual sense that
invariant subsets of $\mathcal{M} \times \Delta^{\mathbb N}$ for $S$ have zero or full
measure.
See \cite{kif88} for equivalent formulations.
The Birkhoff ergodic theorem tells that for an ergodic stationary
measure,
\begin{equation}
\lim_{n \to \infty} \frac{1}{n}\sum_{i=0}^{n-1} \phi
(f^i(x;{\boldsymbol \omega})) = \int_\mathcal{M} \phi(x) d\mu(x)
\end{equation}
for all integrable functions $\phi$ on $\mathcal{M}$ and $\mu\times
\nu^\infty$ almost every point $(x,{\boldsymbol \omega})$. Taking
$\phi=1_A$, the characteristic function of a Borel set $A\subset
\mathcal{M}$, it shows that the relative frequency with which typical orbits
visit $A$ is given by $\mu(A)$.

Write $R^k(\mathcal{M})$ for the space of $C^k$ random diffeomorphisms $f$ on $\mathcal{M}$
(with $f(x;\omega)$ $C^k$ jointly in $x \in \mathcal{M}$ and $\omega\in \Delta$),
 depending on
a random parameter from $\Delta$ through a distribution with a
$C^k$ density function $g$.

Let a random diffeomorphism $f \in R^\infty (\mathcal{M})$ be given.
The existence of finitely many ergodic stationary measures
for $f$
presented in the following theorem,
can be found in \cite[Chapter~5]{doo53}
(valid under more general conditions).
Similar results are contained in \cite{ara00}.
We add statements on the regularity of the stationary measures
valid in our context. Differentiability of stationary densities is also
discussed in \cite{rue90,balyou93}.

\begin{theorem}\label{stationary_measures}
The random diffeomorphism $f \in R^\infty (\mathcal{M})$ possesses
a finite number of ergodic stationary measures
$\mu_1,\ldots, \mu_m$
with mutually disjoint supports $E_i,\ldots, E_m$.
All stationary measures are linear combinations of
$\mu_1,\ldots,\mu_m$.

The support $E_i$ of $\mu_i$ consists of the closure of a finite
number of connected open sets
$C_i^1, \ldots, C_i^p$
that are moved cyclically by $f(\cdot ;\Delta)$.
The density $\phi_i$ of $\mu_i$ is a $C^\infty$ function on $\mathcal{M}$.
\end{theorem}

\noindent {\sc Proof.}
See \cite[Chapter~5]{doo53} for
the existence proof
of cyclically permuted ergodic stationary measures.

Since a point $x$ is mapped to a set $f(x;\Delta)$ diffeomorphic to $\Delta$,
the support $E_i$ is the closure of finitely many connected open sets.
We claim that the closures of the sets $E_i$, $1 \le i \le m$,
are mutually disjoint.
Suppose on the contrary that $\partial E_i$ and $\partial E_j$
with $j \ne i$ have a point $z$ in common.
Then $z$ is mapped by $f(\cdot;\Delta)$ to the injective image
of $\Delta$. By invariance of $E_i$ and $E_j$
under $f(\cdot;\Delta)$,
$f(z;\Delta)$ is contained in both $\bar{E_i}$ and $\bar{E_j}$,
which is not possible.

The regularity statements follow from Proposition~\ref{compact} in Section~\ref{sec_transfer}.
\qed\\

Densities of stationary measures for random endomorphisms
are in general not smooth functions, but are less regular
(see Section~\ref{sec_1D}). Random endomorphisms without
critical points, such as expanding maps, do possess
finitely many smooth stationary densities.
The proofs for random diffeomorphisms extend to cover such random endomorphisms.
The regularity of a stationary density $\phi$ implies that
$\phi$ is flat along the boundary of its support $E$.
By the Birkhoff ergodic theorem (applied to the
characteristic function of a neighborhood of $\partial E$),
this means that typical orbits
are very infrequently found near the boundary of $E$.

We introduce a topology on the space $R^k(\mathcal{M})$
of random diffeomorphisms in order to
be able to compare the dynamics of nearby random diffeomorphisms.
Natural topologies on $R^k(\mathcal{M})$ are
the uniform $C^k$ topologies on
$C^k(\mathcal{M} \times \Delta, \mathcal{M})$.
See e.g. \cite{hir76} for generalities on these topologies.
We will assume $R^k(\mathcal{M})$ to be equipped with this topology.
Note that the alternative approach through discrete Markov processes
suggests a topology using the densities of the transition functions.

Consider $f \in R^\infty (\mathcal{M})$.
Write $\mu_1, \ldots, \mu_m$ for the stationary
measures of $f \in R^\infty(\mathcal{M})$ given by Theorem~\ref{stationary_measures}.

\begin{definition}\label{definition_bifurcation}
A random endomorphism $f \in R^\infty (\mathcal{M})$
is stable if for all $\tilde{f}$ sufficiently close to $f$,
the following two properties are satisfied.
\begin{itemize}
\item
For each $i, 1 \le i \le m$,
the random endomorphism $\tilde{f}$ has a stationary measure
$\tilde{\mu}_i$
whose density  is $C^0$ close to that of $\mu_i$.
\item
The supports of $\tilde{\mu}_i$ and $\mu_i$  are close in the
Hausdorff metric.
\end{itemize}
We speak of a bifurcation, or a bifurcating random endomorphism,
if at least one of these properties is violated.
\end{definition}

\begin{definition}
An ergodic stationary measure  of $f \in R^\infty (\mathcal{M})$
is called isolated or attracting, if there exists an open set
$W$ (an isolating neighborhood) containing
the support $E$ of $\mu$, so that $\overline{f(W;\Delta)} \subset W$ and
$\mu$ is the only ergodic stationary measure of
$f$ with support in $W$.
\end{definition}

For each $\tilde{f}$ close to $f$
is the closure of
$\tilde{f}(W;\Delta)$ contained in $W$. The following stability result shows
that nearby random diffeomorphisms have indeed a unique stationary measure
with support in $W$. For bifurcations where stationary measures vary discontinuously
the condition of being isolated must therefore be violated.
The proof of the following theorem is found in Section~\ref{sec_stable}.

\begin{theorem}\label{isolated}
Let $\mu$ be an isolated ergodic stationary measure of $f \in R^\infty (\mathcal{M})$
with density $\phi$ with isolating neighborhood $W$.
Then each $\tilde{f} \in R^\infty (\mathcal{M})$ sufficiently close to
$f$ possesses a unique ergodic stationary measure $\tilde{\mu}$
with support
in $W$. The density $\tilde{\phi}$ of $\tilde{\mu}$ is $C^\infty$ close
to $\phi$.
\end{theorem}

Note though that the above theorem leaves open the possibility that
the supports $\tilde{E}$ of $\tilde{\mu}$ and $E$ of $\mu$
are not close
in the Hausdorff metric.
An illustrative example
of this phenomenon is described in Section~\ref{sec_circle}.

The following theorem
establishes that stable random diffeomorphisms
are generic.
Its proof is in Section~\ref{sec_stable}.
The argument also shows that random
diffeomorphisms with
a locally constant number of smoothly varying stationary densities
(ignoring variations in their support)
form an open and dense subset of $R^\infty (\mathcal{M})$.

\begin{theorem}\label{stable}
The set of stable random diffeomorphisms in $R^\infty (\mathcal{M})$
contains a countable intersection of open and dense sets.
\end{theorem}

A thorough description of the dynamics of random circle diffeomorphisms
is in Section~\ref{sec_circle}.
Random diffeomorphisms, and even random endomorphisms, on the circle
are shown to form an open and dense set, see Theorem~\ref{stable_1D}.

\subsection{Families of random diffeomorphisms}\label{subsec_families}

Bifurcations are best studied in families depending on finitely many parameters.
We will consider families of random diffeomorphisms depending on a single real
parameter,
where we have the goal to focus
on bifurcations that typically occur varying one parameter.

\begin{definition}
A {\em smooth family of random endomorphisms}
is a family of  random endomorphisms $\{f_a\}$
depending on parameters $a$,
so that $f_a(x;\omega)$ depends smoothly on $(x,\omega,a)$.
A {\em smooth family of random diffeomorphisms} is
a smooth family of random endomorphisms
where each map $f_a(\cdot;\omega)$ is a diffeomorphism.
\end{definition}

\begin{remark}
Alternatively, one can explicitly include
noise densities $g_a$ in the definition (considering pairs
$(f_a,g_a)$) with $g_a(\omega)$ varying smoothly with $(\omega,a)$.
Compare Remark~\ref{remark_g}.
For convenience we consider fixed noise densities, but
completely analogous results hold if the noise densities
are allowed to vary with $a$.
\end{remark}

Consider a smooth one parameter family $\{ f_a \}$ of random diffeomorphisms,
with $a$ from an interval $I$.
Consider a parameter value $a_0$ and an ergodic stationary measure $\mu_{a_0}$
with support $E_{a_0}$.
The following result extends Theorem~\ref{isolated}, providing
an analogous statement in the context of families.
If $\mu_{a_0}$ is an isolated  ergodic stationary measure then there are
ergodic invariant measures $\mu_a$ for $a$ near $a_0$ with nearby densities.

\begin{theorem}\label{continuation}
Suppose $\mu_{a_0}$
is an isolated ergodic stationary measure.
Then the stationary density $x\mapsto \phi_a(x)$ of $\mu_{a}$
depends $C^\infty$ on $(x,a)$.
\end{theorem}

See Section~\ref{sec_parameter} for the proof.
We stress again that the support $E_a$ of $\mu_a$ can still
vary discontinuously in the Hausdorff metric with $a$.
The number of components of the support of the stationary measure
can also change, while the stationary density varies smoothly.

Consider a smooth function $\phi$ with support on an isolating
neighborhood $W$ for $\mu_{a_0}$ and compute averages of $\phi$
along orbits $f^k_a(x;{\boldsymbol \omega})$. By the Birkhoff
ergodic theorem, for typical initial points $x \in W$ and noise
sequences ${\boldsymbol \omega}$, the averages lie on a smooth
function of $a$ for $a$ near $a_0$.

The two types of bifurcation distinguished in Definition~\ref{definition_bifurcation}
gives rise to a particular dynamical phenomenon
associated to either intermittency or transients.
Consider a family $\{f_a\}$ of random diffeomorphisms in $R^\infty(\mathcal{M})$,
with $a$ from an open interval $I$.
Suppose that $a_0\in I$ is a bifurcation value for $\{f_a\}$ involving
a stationary measure $\mu$. Write $\phi$ for the density of $\mu$.
Analogies with deterministic dynamics suggest the following two definitions.
In Section~\ref{sec_1D} we will see that in typical one parameter
families of random interval or circle endomorphisms bifurcations
are isolated and of these two types.

\begin{definition}
The bifurcation at $a_0$ is called an {\em intermittency bifurcation} if
there is a stationary density $\phi_a$ for $\{f_a\}$ with $\phi_{a_0} = \phi$ and
depending continuously on $a$, so that the support $E_a$ of $\phi_a$ varies with $a$,
for $a$ near $a_0$, as follows.
\begin{itemize}

\item $E_a$ varies continuously for $a$ from one side of $a_0$. Without loss
of generality, we assume this to be the case for $a<a_0$.

\item $E_a$ is discontinuous at $a=a_0$ and $E_a$ contains an open set disjoint from
$E_{a_0}$ for $a>a_0$.

\end{itemize}
\end{definition}

An orbit piece outside a small neighborhood $W$ of $E_{a_0}$ is called a burst.
Out of the substantial literature on intermittency in dynamics, we point to references
\cite{pomman80,eckthowit81,hirhubsca82,greottromyor87,homyou02,homyou05}.

\begin{definition}
The bifurcation at $a_0$ is called a {\em transient bifurcation} if
there is a stationary density $\phi_a$ for $\{f_a\}$ with $\phi_{a_0} = \phi$
for $a$ close to $a_0$ from one side of $a_0$
(without loss
of generality, we assume this to be the case for $a < a_0$),
so that
\begin{itemize}

\item $\phi_a$ and its support $E_a$ vary continuously with $a$, for $a \le a_0$.

\item there is no stationary density near $\phi_{a_0}$ for $a$ close to $a_0$ and $a > a_0$.

\end{itemize}
\end{definition}

\begin{figure}[htb]
\centerline{\hbox{
\epsfig{figure=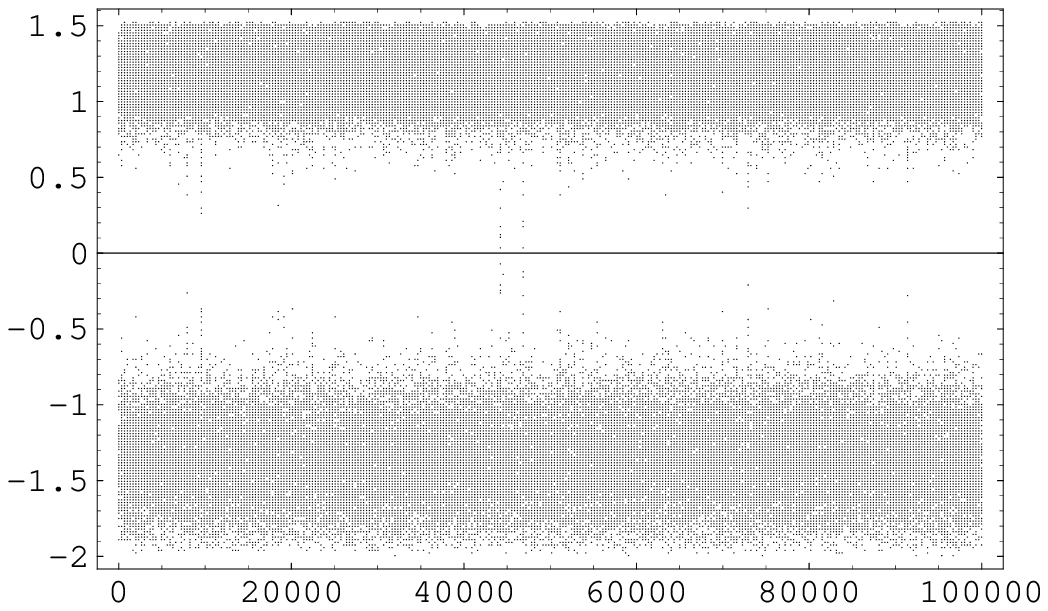,height=5cm,width=7cm} \hspace{0.5cm}
\epsfig{figure=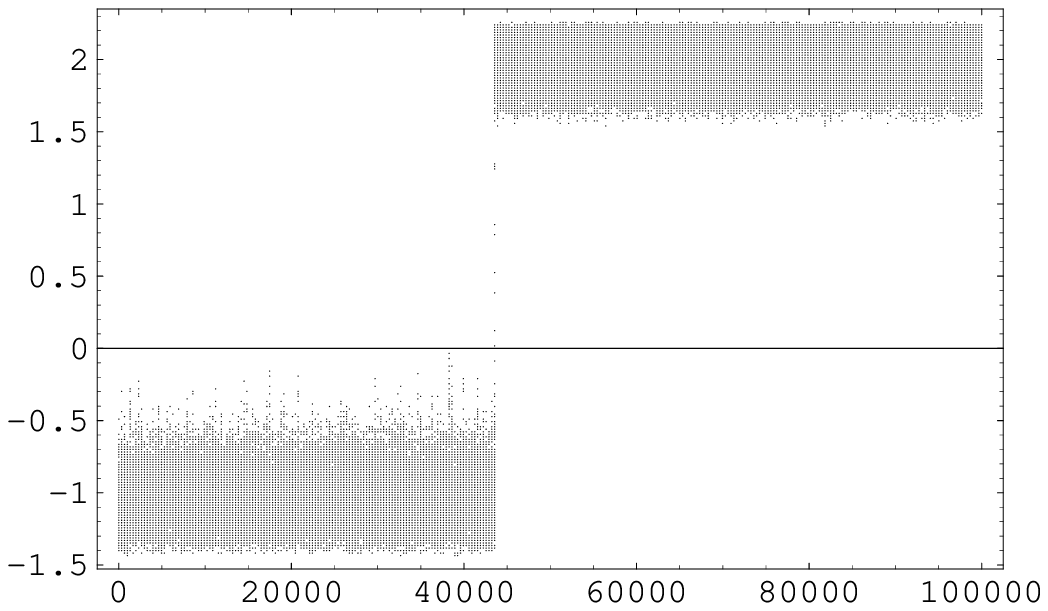,height=5cm,width=7cm}
}}
\caption{\small
Typical times series for intermittent dynamics on the left and
transient dynamics on the right.
The time series are for bifurcation values after the bifurcation
took place.
The intermittency bifurcation involves interval diffeomorphisms
with a single stationary measure;
the support consisting of two intervals for
a period two cycle bifurcates to form a single interval.
In the transient bifurcation one stationary measure out of two
stationary measures existing previous to the bifurcation disappears.
\label{fig_timeseries}}
\end{figure}

We end this section with a quantitative
description of time series near intermittent or transient random
bifurcations.
We do this through the estimation of the expected escape time
from a small neighborhood of the support
of the bifurcation stationary measure.
Estimating the speed of decay of correlations as function of a parameter
gives further details of changes through bifurcations.

As before, let $\{f_a\}$ be a family of random diffeomorphisms in
$R^\infty(\mathcal{M})$, with $a$ from an open interval $I$. Let $\mu_{a_0}$
be a stationary measure of $f_{a_0}$ for some $a_0 \in I$ and let
$W$ be an open neighborhood of the support $E_{a_0}$ of
$\mu_{a_0}$ such that no other stationary measure has support
intersecting $\overline{W}$. If $f_{a_0}$ is stable, then there is
a unique stationary measure $\mu_a$ that is the continuation of
$\mu_{a_0}$. The stationary measure $\mu_a$ has its support in $W$
for $a$ near $a_0$. If $a_0$ is a bifurcation value for $\{f_a\}$,
iterates $f^k_a(x;{\boldsymbol \omega})$, for certain $x \in W$
and $a$ near $a_0$, may leave $W$.

For $x \in W$ and ${\boldsymbol \omega} \in \Delta^{\mathbb N}$,
define
\begin{equation}\label{chi}
\chi_a(x,{\boldsymbol \omega}) = \min \{k \;|\; f^k_a
(x;{\boldsymbol \omega}) \not \in W).
\end{equation}
The following result shows how the average escape
time from a neighborhood of the support of a bifurcating
stationary measure
is 
more than polynomially large in an
unfolding parameter. This makes it difficult to accurately
establish the bifurcation parameter value
using finite data, even in numerical simulations.
It explains the occurrence of very long transients
near a transient bifurcation and the very irregular
occurrence of bursts in intermittent time series.
The proof, in Section~\ref{sec_escapeII}, relies on
the construction of conditionally stationary measures in Section~\ref{sec_escape}.

\begin{theorem}\label{theorem_exittimes_smooth}
For each $k>0$ there is a constant $C_k>0$ so that
$$\int_{W} \int_{\Delta^{\mathbb N}}
       \chi_a (x ,{\boldsymbol \omega}) d\nu^\infty({\boldsymbol \omega}) dm(x) \ge C_k
\left| a-a_0 \right|^{-k}.
$$
\end{theorem}

A single random map with an isolated measure supported on a single component
has exponential decay of correlations as precised in the following proposition.
The interest from our perspective
in computing the speed of decay of correlations lies
in the study of bifurcations where the support
of a stationary measure has several components merging.
This will be discussed below.
Proofs of the following statements are in Section~\ref{sec_correlations}.
The reader can consult \cite{via97,bal00} for background on decay of correlations.
Write
$$
U^n \psi (x)  =  \int_{\Delta^n} \psi \circ f^n (x;\omega_1,\ldots,\omega_n)
                        d\nu(\omega_1)\cdots d\nu(\omega_n).
$$

\begin{proposition}\label{prop_corr}
Let $f$ be a random map with an isolated stationary measure $\mu$
with connected support. Let $W$ be an isolating neighborhood for $\mu$.
Take $\varphi,\psi \in {\cal L}^2 (W)$.
Then
$$
\left|
\int_\mathcal{M} \varphi(x) U^n \psi(x) dm(x) - \int_\mathcal{M} \varphi(x) dm(x) \int_\mathcal{M} \psi(x) d\mu(x)
\right|
\le C \eta^n
$$
for some $C>0$, $0<\eta<1$.
\end{proposition}

Note that the exponential decay of correlations holds for
observables $\varphi,\psi \in {\cal L}^2(W)$, thus including characteristic
functions of open sets. On the other hand, our definition involves
an average over noise sequences. The following remark addresses this point.

\begin{remark}\label{rem_decay}
For a fixed noise sequence ${\boldsymbol \omega} \in
\Delta^{\mathbb N}$ and nonnegative observables $\varphi,\psi \in {\cal
L}^2(W)$, consider
\begin{equation}\label{corr_fixedomega}
\left|
\int_\mathcal{M} \varphi(x) \psi\circ f^n (x;\omega_1,\ldots,\omega_n)(x) dm(x) - \int_\mathcal{M} \varphi(x) dm(x) \int_\mathcal{M} \psi(x) d\mu(x)
\right|.
\end{equation}
By  Proposition~\ref{prop_corr}, the integral over $\Delta^n$ of this expression is bounded by $C \eta^n$.
Necessarily,
(\ref{corr_fixedomega}) is exponentially small in $n$ for
$(\omega_1,\ldots,\omega_n)$ outside an exponentially small set.
Indeed, choose $\tilde{\eta}$ slightly larger than $\eta$.
Then (\ref{corr_fixedomega}) is larger than $\tilde{\eta}^n$ only on a 
set $S_n\subset \Delta^n$ with $\nu^n (S_n) < C \eta^n / \tilde{\eta}^n$. 
\end{remark}

We return to a family $\{f_a\}$, $a \in I$, of random maps.
Assume that $f_a$ has an isolated measure $\mu_a$ for all $a \in I$
with an isolating neighborhood $W$.
Suppose $a_0 \in I$ is a bifurcation value for an intermittency bifurcation
so that
\begin{itemize}
\item
the support of $\mu_a$ consists of $k$ components for $a \le a_0$,
\item
the support of $\mu_a$ consists of a single component for $a >a_0$.
\end{itemize}
We incorporate the dependence of $U^n$ on $a$ into  the notation by writing $U^n_a$.

\begin{theorem}\label{theorem_corr}
Let $\{f_a\}$ be as above.
Take $\varphi,\psi \in {\cal L}^2(W)$.
There are a constant $C>0$ and a smooth function $a \mapsto \eta_a$
with $\eta_a = 1$ for $a<a_0$ and $\eta_a < 1$ for $a > a_0$, so that
$$
\left|
\int_\mathcal{M} \varphi(x) U^n_a \psi(x) dm(x) - \int_\mathcal{M} \varphi(x) dm(x) \int_\mathcal{M} \psi(x) d\mu_a(x)
\right|
\le C \eta_a^n
$$
for $a > a_0$.
\end{theorem}

The smoothness properties of $\eta_a$ imply that
$\eta_a - 1$ is a flat function of $a$ at $a=a_0$.
As in the discussion of escape times, that shows how slowly
the bifurcation manifests itself in time series
when moving the parameter $a$.
Remark~\ref{rem_decay} also applies in the parameter dependent context.


\section{Transfer operators}\label{sec_transfer}

Let $f$ be a random diffeomorphism on the manifold $\mathcal{M}$.
Associated to $f$ is the stochastic transition function,
$$
P(x,A) = \int_{\{\omega \;|\; f(x;\omega) \in A\}}  d\nu(\omega) =
\int_\Delta  1_A (f(x;\omega))d\nu(\omega),
$$
for Borel sets $A\subset \mathcal{M}$.
Write $h_x(\omega) = f(x;\omega)$, note that $h$ maps $\Delta$ injectively
onto $U_x = f(x;\Delta)$.
The density
\begin{equation}\label{kernel}
k(x,y) = d \left( h_x \right)_* \nu / dm
\end{equation}
of $P(x,\cdot)$
vanishes for $y$ outside $U_x$ and is a smooth function on its support.

Write ${\cal L}^1(\mathcal{M})$ for the space of integrable functions on $\mathcal{M}$.
Define the transfer operator $L$
acting on ${\cal L}^1(\mathcal{M})$ by
\begin{equation}\label{transferA}
\int_A L\phi (x) dm(x) = \int_\mathcal{M} P(x,A) \phi(x) dm(x).
\end{equation}
The transfer operator is a positive linear operator.
A stationary density is a fixed point of $L$.

Define 
\begin{equation}
V_x = \{ z \in \mathcal{M} \;|\; x \in f(z;\Delta) \},
\end{equation}
which is the set
of points in $\mathcal{M}$ that are mapped to $x$ by some random map.
The assumptions on the random diffeomorphism $f$
imply that the equation $f(x;\omega) = y$ can be solved for $x$ as a diffeomorphic
map of $\omega$. Therefore
$V_x$ is diffeomorphic to $\Delta$ and thus
a domain with piecewise smooth boundary, depending
smoothly on $x$.

Write $P_{f(\cdot;\omega)}$ for the Perron-Frobenius operator, defined by
\begin{equation}\label{P-F}
\int_A P_{f(\cdot;\omega)} \phi(x) dm(x) = \int_{f^{-1} (A;\omega)} \phi(x) dm(x)
\end{equation}
for Borel sets $A$.
That is, for
a measure $\mu$ with density $\phi = d \mu / dm$,
$f_* \mu$ has density $P_f \phi = d f_* \mu / dm$ (see e.g. \cite{lasmac94}).
The following lemma gives the transfer operator as an average
over the random parameters $\omega$ of
the Perron-Frobenius operators for $f(\cdot;\omega)$ and gives an
equivalent formulation as an integral over the state space $\mathcal{M}$.

\begin{lemma}\label{lemma_transfer}
The transfer operator is given by
\begin{equation}\label{transfer_average}
L \phi(x) = \int_\Delta P_{f(\cdot;\omega)} \phi (x)  d\nu(\omega),
\end{equation}
or,
\begin{equation}\label{transfer}
L \phi (x) = \int_{V_x} k(y,x) \phi(y) dm(y)
\end{equation}
\end{lemma}

\noindent {\sc Proof.}
By (\ref{transferA}), for a continuous function $\psi$,
\begin{eqnarray*}
\int_\mathcal{M} \psi(x) L\phi (x) dm(x) & = &
\int_\mathcal{M} \int_{\Delta} \psi(f (x;\omega)) \phi (x) d \nu(\omega) dm(x)
\end{eqnarray*}
Calculate
\begin{eqnarray*}
\int_\mathcal{M} \int_{\Delta} \psi(f (x;\omega)) \phi (x) d \nu(\omega) dm(x)
& = &
\int_{\Delta} \int_\mathcal{M} \psi(f (x;\omega)) \phi (x) dm(x) d \nu(\omega)
\\
& = &
\int_{\Delta} \int_\mathcal{M} \psi(x) P_{f(\cdot;\omega)} \phi(x) dm(x) d \nu(\omega)
\end{eqnarray*}
This implies (\ref{transfer_average}).
Alternatively,
\begin{eqnarray*}
\int_\mathcal{M} \int_{\Delta} \psi(f (x;\omega)) \phi (x) d \nu(\omega) dm(x)
& = &
\int_\mathcal{M} \int_{U_x} \psi(y) \phi (x) k(y,x) dm(y) dm(x)
\\
& = & 
\int_\mathcal{M}  \int_{V_y} \psi(y) \phi (x) k(y,x) dm(x) dm(y)
\end{eqnarray*}
proves (\ref{transfer}).
\qed

\begin{remark}
The transfer operator $L$ preserves integrals, as the following
computation shows.
\begin{eqnarray*}
\int_\mathcal{M} L \phi (y) dm(y) & = & \int_\mathcal{M} \int_{V_y} k(x,y) \phi(x) dm(x) dm(y)
\\
& = &  \int_\mathcal{M} \int_{U_x} k(x,y) \phi(x) dm(y) dm(x)
\\
& = & \int_\mathcal{M} \phi(x) dm(x),
\end{eqnarray*}
because $\int_{U_x} k(x,y) dm(x) = 1$.
\end{remark}

Iterating $L$ gives
\begin{eqnarray*}
L^2 \phi (x) & = & \int_{V_x} k(z,x) L\phi(z) dz 
\\ & = &  
\int_{V_x} \int_{V_z}  k(z,x) k(y,z)  \phi(y) dy dz 
\\
& = & 
\int_{f^{-1}(V_x,\Delta)} \int_{V_x}  k(z,x) k(y,z) dz \phi(y) dy
\end{eqnarray*}
which is of a similar form,
namely $\int_{f^{-1}(V_x,\Delta)} k_2(y,x) \phi(y) dy$ with
$k_2(y,x) = \int_{V_x}  k(z,x) k(y,z) dz$,
as $L \phi (x)$. Inductively similar expressions are 
derived for higher iterates of $L$.

Denote by
${\cal D}(\mathcal{M}) = \{ \phi\in {\cal L}^1(\mathcal{M}) \;|\; \phi \ge 0,  \int_\mathcal{M} \phi(x) dm(x) = 1\}$
the space of densities on $\mathcal{M}$.
The above remark shows that $L$ maps ${\cal D}$ into itself.
Smoothness of stationary densities is obtained by showing that $L$ maps
a space of smooth densities into itself.

\begin{proposition}\label{compact}
The transfer operator $L$ maps $C^k(\mathcal{M})$ into itself and
is a compact operator on  $C^k(\mathcal{M})$.

The number 1 is an eigenvalue of $L$
with equal algebraic and geometric multiplicity $m \ge 1$.
The densities $\phi_1,\ldots,\phi_m$ provide a basis
of eigenfunctions with mutually disjoint support.
Each eigenfunction $\phi_i$ is $C^\infty$ and its support
consists of a finite number $c_i$ of connected components.
\end{proposition}

\noindent {\sc Proof.}
Theorem~\ref{stationary_measures}
gives the $m$ invariant densities $\phi_1, \ldots,\phi_m$.
The geometric multiplicity of the eigenvalues 1 is equal to $m$.
Since $L$ preserves the ${\cal L}^1$ norm,
the algebraic and geometric multiplicity of 1 are equal.
To see this, suppose on the contrary that there is a nontrivial vector in
$\ker (L - I)^2 \backslash \ker(L-I)$.
Elementary linear algebra gives the existence of a sequence of vectors
$\psi_n$ inside $\ker (L-I)^2$ converging to an eigenvector $\phi$, such that
$\lim_{n\to \infty} L^n \psi_n = 2 \phi$. Indeed, take $\psi \in \ker (L-I)^2$ with
$L \psi = \phi + \psi$ and let $\psi_n = \frac{1}{n}\psi$. From
$L^n(\psi) = n \phi + \psi$ it follows that
$L^n (\phi + \psi_n) =  2 \phi  + \psi_n$, which converges
to $2 \phi$ if $n\to \infty$.
This contradicts
the preservation of the ${\cal L}^1$ norm by $L$.

There can be no additional eigenvectors of $L$ that do not
correspond to linear combinations of densities.
Namely, suppose $\phi$ is an eigenvector taking both positive and
negative values. Write $\phi = \phi_+ - \phi_-$ for nonnegative functions
$\phi_+ = \max\{\phi,0\}$, $\phi_- = \max\{-\phi,0\}$.
If $\phi$ is not a  linear combinations of densities,
the supports of $\phi_+, \phi_-$ cannot be invariant.
Since $L$ is positive and preserves the ${\cal L}^1$ norm,
$(L \phi)_+ < L \phi_+$ so that $\phi$ cannot be an eigenvector.

We will show that $L$ maps ${\cal L}^1 (\mathcal{M})$ into $C^0(\mathcal{M})$ and
$C^i(\mathcal{M})$ into $C^{i+1}(\mathcal{M})$. From this it follows that $\phi_i \in C^{\infty} (\mathcal{M})$.
Take $\psi \in {\cal L}^1 (\mathcal{M})$.
Use a chart to identify a neighborhood of $x \in \mathcal{M}$ with an open set in ${\mathbb R}^n$.
With $h$ a small vector in ${\mathbb R}^n$, consider
\begin{eqnarray}
\nonumber
 L\psi (x+h) - L\psi (x)  & = &
\int_{V_{x+h}} k(y,x) \psi (y) dm(y) -  \int_{V_{x}} k(y,x) \psi (y) dm(y)
\\
\nonumber
& = &
 \int_{V_{x+h}\cap V_{x}} ( k(y,x+h) - k(y,x) ) \psi (y) dm(y)
\\
\nonumber
& &
+
 \int_{V_{x+h} \backslash (V_{x+h}\cap V_{x})} k(y,x+h) \psi(y) dm(y)
\\
\label{continuityL}
& &
-
\int_{V_{x} \backslash (V_{x+h}\cap V_{x})} k(y,x) \psi(y) dm(y).
\end{eqnarray}
The first term on the right hand side is small for $h$ small
by continuity of $k$ and integrability of $\psi$.
The other two terms are small for $h$ small
by the continuous dependence of $V_x$ on $x$.
Continuity of $L\psi$ follows.
Suppose next that $\psi \in C^0(\mathcal{M})$ and consider
$\frac{1}{|h|} \left( L\psi (x+h) - L\psi (x)  \right)$
This equals the right hand side of (\ref{continuityL}) divided by $|h|$.
Note that
$$
\lim_{h\to 0} \frac{1}{|h|}  \int_{V_{x+h}\cap V_{x}} ( k(y,x+h) - k(y,x) ) \psi (y) dm(y)
= \int_{V_{x}} \frac{\partial}{\partial x} k(y,x) \frac{h}{|h|} \psi(y) dm(y)
$$
is a continuous function of $x$.
To check continuity of the remaining two terms it suffices to
do a local calculation by covering the boundary of $V_x$ by finitely many balls
and using a partition of unity.
Without loss of generality we may assume that
near a smooth part of the boundary of $V_x$,
$V_x$ is bounded from below by the graph of a continuously differentiable
function $H_x : [0,1]^{n-1}\to {\mathbb R}$.
We may also assume that $m$ equals Lebesgue measure on ${\mathbb R}^n$.
Then
\begin{eqnarray*}
\lim_{h\to 0} \frac{1}{|h|}
\int_{V_{x} \backslash (V_{x+h}\cap V_{x})} k(y,x) \psi(y) dy & = &
\lim_{h \to 0} \frac{1}{|h|}
\int_{[0,1]^{n-1}}  \int_{H_x(y_1)}^{H_{x+h}(y_1)} k(y_1,y_2,x)  \psi(y) dy_2 dy_1
\\
& = &
\int_{[0,1]^{n-1}} \frac{\partial}{\partial x} H_x (y_1) \frac{h}{|h|}
k(y_1,H_x(y_1),x) \psi(y_1,H_x(y_1)) d y_1
\end{eqnarray*}
is a continuous function of $x$.
The contribution near the
finitely many points where $V_x$ is not smooth
vanishes in the limit $h \to 0$.
Summarizing, $D (L\psi)$ has an expression of the form
\begin{equation}\label{derivative}
D(L\psi)(x) = \int_{V_x} \frac{\partial}{\partial x} k(y,x) \psi(y) dm(y)
 + \int_{\partial V_x} n(y,x) k(y,x) \psi(y) dS(y)
\end{equation}
where $n(y,x)$ measures the change of $\partial V_x$ in the direction of the
unit normal vector to $V_x$
and $dS$ is the volume on $\partial V_x$.
We remark that 
the formula is a variant of the transport theorem 7.1.12 and the Gauss theorem 7.2.9 in
\cite{abrmarrat83}.
It follows that $L\psi$ is continuously differentiable if $\psi$ is continuous.
Higher order derivatives are computed inductively. 
This gives that $L\psi \in C^{k+1}(\mathcal{M})$ for $\psi \in C^k(\mathcal{M})$.

We prove compactness on $C^k(\mathcal{M})$ by
modifying the argument in \cite{mac91}.
Let $B^k(\mathcal{M})$ be the unit sphere in $C^k(\mathcal{M})$.
Consider first $k=0$.
By the Arzela-Ascoli theorem, compactness of $L$ on $C^0(\mathcal{M})$ follows
from the following two properties (compare \cite{zee88}),
\begin{itemize}
\item for all $x \in \mathcal{M}$, $\{ |L\psi(x)| \;|\; \psi \in B^0(\mathcal{M}) \}$ is bounded,
\item $LF$ is equicontinuous.
\end{itemize}
For $\psi \in B^0(\mathcal{M})$, $|L \psi (x)| \le \int_\mathcal{M} k(y,x) dm(y)$. This is a continuous
function of $x$ and hence bounded. This proves the first item.
The above computations showing that $L\psi$ is continuously differentiable
also show that
$\| D (L\psi) (x)\|$ is uniformly bounded on $B^0(\mathcal{M})$.
This proves that $LF$ is equicontinuous.
Compactness in $C^k(\mathcal{M})$ follows similarly by noting that
\begin{itemize}
\item for all $x \in \mathcal{M}$, $i\le k$,  $\{ \|D^i (L\psi)(x)\| \;|\; \psi \in B^k(\mathcal{M}) \}$ is bounded,
\item $\| D^{k+1} (L\psi) (x)\|$ is uniformly bounded on $B^k(\mathcal{M})$.
\end{itemize}
\qed

\begin{remark}\label{remark_l2}
The transfer operator $L$ is compact
on the space ${\cal L}^2(\mathcal{M})$ of quadratic
integrable functions on $\mathcal{M}$ \cite[Section~X.2]{yos80} and \cite{delljun99}.
The proof of Proposition~\ref{compact},
demonstrating that the transfer operator
increases regularity of functions,
implies that the spectrum of $L$ on
${\cal L}^2(\mathcal{M})$ equals that of $L$ on $C^k(\mathcal{M})$.
\end{remark}

\begin{remark}\label{remark_tobias}
The spectral radius of $L$ is 1, the eigenvalue 1 occurs with multiplicity
$m$ equal to the number of stationary measures. The peripheral spectrum
on the unit circle consists of eigenvalues $e^{2 \pi i/p}$, 
$ 0 \le i < p$, for each $p$
occurring as the number of connected components of a stationary
measure.  See \cite{gay01} and \cite[Theorem~V.4.9]{sch74}
for a proof.
\end{remark}

\begin{proposition}\label{Lcontinuous}
The transfer operator $L$ as a linear map on $C^k(\mathcal{M})$ or ${\cal L}^2(\mathcal{M})$
depends continuously on $f \in R^k(\mathcal{M})$.
\end{proposition}

\noindent{\sc Proof.}
Consider $\tilde{f}$ near $f$. Write $\tilde{L}$ and $L$ for the
corresponding transfer operators.
We need to prove that $\tilde{L} - L$ has small norm.
Consider the transfer operators operating on $C^k(\mathcal{M})$
(continuity on ${\cal L}^2(\mathcal{M})$ is treated analogously).
The transfer operator $\tilde{L}$ is given as
$\tilde{L}\phi (x) = \int_{\tilde{V}_x} \tilde{k} (y,x) \phi(y) dm(y)$.
For $\phi \in B^k(\mathcal{M})$, the unit sphere in $C^k(\mathcal{M})$,
$$\tilde{L} \phi (x) - L\phi (x) =
\int_{\tilde{V}_x} \tilde{k}(y,x) \phi(y) dm(y) - \int_{V_x} k(y,x) \phi(y) dm(y)$$
is small, uniformly in $x$,
since $\tilde{k}(y,x)$ is close to $k(y,x)$  on $\tilde{V}_x \cap V_x$
and $\tilde{V}_x$ is close to $V_x$.
The derivative $D (L\phi)$ is given by (\ref{derivative}). An analogous
formula holds for $D (\tilde{L}\phi)$.
Since the functions and sets involved in the two formulas
for $D(\tilde{L}\phi)$ and $D(L\phi)$ are close,
$D(\tilde{L} \phi) (x)$ is uniformly close to $D(L\phi) (x)$.
Closeness of higher order derivatives, up to order $k$, is treated analogously.
Continuity on ${\cal L}^2(\mathcal{M})$ is proved analogously.
\qed

\begin{remark}\label{remark_kato}
Consider two nearby random diffeomorphisms $f$ and $\tilde{f}$
from $R^k(\mathcal{M})$. Write $L$ and $\tilde{L}$ for the corresponding transfer operators
on $C^k(\mathcal{M})$. Let  $\lambda_1, \ldots, \lambda_l$ be a finite set of eigenvalues
for $L$ and denote by $F$ the sum of the corresponding generalized eigenspaces.
Then $\tilde{L}$ possesses a
nearby set of eigenvalues $\tilde{\lambda}_1, \ldots, \tilde{\lambda}_l$.
The sum $\tilde{F}$ of the corresponding
generalized eigenspaces is a small perturbation of $F$ (in the sense that $F$ and $\tilde{F}$
have nearby bases). See \cite[Theorem~IV.3.16]{kat66}.
\end{remark}


\section{Stable random diffeomorphisms}\label{sec_stable}

This section contains the proofs of Theorem~\ref{isolated} on
stability of isolated stationary measures and
Theorem~\ref{stable} establishing generic stability
of random diffeomorphisms. 

For the restriction of random maps to an isolating neighborhood
$W$ we consider the transfer operator acting on functions
vanishing outside $W$ and at the boundary of $W$.
Write  
$$
C^k_0(W) = \{ f \in C^k(\mathcal{M})\;|\; \text{the support of }f\text{ is contained in }\bar{W}\}.
$$ 
Then $L$ acting on $C^k_0(W)$ is well defined.
The results in the previous section hold for $L$ acting on $C^k_0(W)$. \\

\noindent {\sc Proof of Theorem~\ref{isolated}.}
Recall that the closure of
$f(W;\Delta)$ is contained in the isolating neighborhood $W$.
This property extends to
random diffeomorphisms sufficiently close to $f$.
Restrict the map $x \mapsto f(x;\omega)$ to $W$
and
and consider
the transfer operator $L$ acting on $C^k_0(W)$.
Then $L$ has a single eigenvalue 1.
Since the spectrum of the transfer operator varies continuously with the
random diffeomorphism at $f$, the transfer operator
corresponding to each nearby random diffeomorphism
possesses a single eigenvalue 1.
The corresponding eigenvector is near $\phi$. \qed 

\begin{lemma}\label{lem_stable}
Write $L_f$ for the transfer operator on $C^k(\mathcal{M})$ for
$f \in R^\infty(\mathcal{M})$.
Densities of stationary measures vary continuously with $f \in R^\infty(\mathcal{M})$ 
at a random diffeomorphism $\bar{f}$ precisely if 
the multiplicity of the eigenvalue 1 for $L_f$ is locally constant in $f$ for $f$ near
$\bar{f}$.
\end{lemma}

\noindent {\sc Proof.}
Consider $\bar{f} \in R^\infty(\mathcal{M})$ with an eigenvalue 1 of multiplicity $m$.
Let $\bar{\mu}_1 , \ldots, \bar{\mu}_m$ be the ergodic stationary measures
with densities $\bar{\phi}_1 , \ldots , \bar{\phi}_m$.
Write $F_{\bar{f}}$ be the direct sum
of the lines spanned by
$\bar{\phi}_1 , \ldots , \bar{\phi}_m$.
By Remark~\ref{remark_kato}, 
the transfer operator for any
 $f \in R^\infty(\mathcal{M})$ sufficiently close to $\bar{f}$
possesses a $m$-dimensional invariant
linear space $F_f$ that is the continuation of $F_{\bar{f}}$.

The spectrum of $L_{f}$ restricted to $F_{f}$
is in general close to 1.
Suppose now that all eigenvalues equal 1.
Then  $\phi \in  F_{f}$ implies $L_{f} \phi = \phi$.
Write $\phi = \phi^+ - \phi^-$ with $\phi^+ = \max \{0,\phi\}$ and $\phi^- = \max \{0,-\phi\}$
the positive and negative parts of $\phi$.
Because $L_{f} \ge 0$ and $L_{f}$ preserves the ${\cal L}^1$ norm,
$(L_{f} \phi)^+ < L_{f} \phi^+$ precisely if $\phi^-$ and $\phi^+$ are not invariant.
Thus $\phi^+$ and $\phi^-$ are necessarily invariant.
It follows that invariant densities are obtained by
taking positive parts of invariant eigenfunctions.
This way $m$ invariant densities for $f$ near those of $\tilde{f}$
can be obtained, proving the lemma.
\qed \\

\noindent {\sc Proof of Theorem~\ref{stable}.}
Consider diffeomorphisms on an open neighborhood $U$ of
a random diffeomorphism $\bar{f} \in R^\infty (\mathcal{M})$.
Write $m$ for the multiplicity of the eigenvalue 1
for the transfer operator corresponding to $\bar{f}$.
There is a neighborhood $D$ of 1 in the complex plane, so that
for $U$ small enough, each $f \in U$ has $m$ eigenvalues counting
multiplicity in $D$.
Let $F$ denote the $m$ dimensional invariant linear space
corresponding to these eigenvalues.
Consider the map that assigns to $f \in U$
the union of the support of all functions in $F$.
By the continuous dependence of $F$ on $f$, 
this is a lower semicontinuous set valued mapping
and therefore continuous on a set $B_2 \subset U$
of Baire second category \cite{for51}.

Consider the map that assigns to random diffeomorphisms
$f \in R^\infty (\mathcal{M})$
the multiplicity $m(f)$ of the eigenvalue 1
for the corresponding transfer map.
By the continuous dependence of eigenvalues of the transfer map on $f$, 
the map $m$ is upper semicontinuous. Since $m$ takes on finitely many values,
it is continuous on an open and dense subset of $R^\infty (\mathcal{M})$.
Indeed, consider $A_n = \{ f \in R^\infty (\mathcal{M}) \;|\; m(f) < n\}$.
The set of points of continuity of $m$,
in the vicinity of some map in $R^\infty (\mathcal{M})$,
equals the intersection of a finite collection of
open and dense sets $A_n \cup \left( R^\infty (\mathcal{M}) \backslash \overline{A_n} \right)$, namely
with $n$ ranging over a finite set of positive integers.

With reference to Lemma~\ref{lem_stable}, the two above items combined prove the theorem.
\qed

\section{Auxiliary parameters}\label{sec_parameter}

Consider a smooth one parameter family of random diffeomorphisms
$x\mapsto f_a (x;\omega)$ depending on $a$ from an open interval $I$ in ${\mathbb R}$.
The transition map $P$ and its density $k$ depend on $a$, we write $P_a$ and $k_a$.
The support $U_{x,a}$ of $k_a$ is assumed to vary smoothly with $x$ and $a$.
The density $k_a(x,y)$ is a smooth function of
$(a,x,y) \in \cup_{a,x} \{a\} \times \{x\} \times U_{x,a}$ in the sense that it
can be extended to a smooth function on an open neighborhood.
Let $L_a$ denote the transfer operator for $f_a$, given by
\begin{equation}\label{Lwitha}
L_a\phi(x) = \int_{V_{x,a}} k_a(y,x) \phi(y) dm(y).
\end{equation}
The domain of integration
$V_{x,a} = \{ y \in \mathcal{M} \;|\; y \in f_a(x;\Delta) \}$
depends smoothly on $(x,a)$.

%

\begin{proposition}\label{prop_smoothnessLa}
For $r \ge 0$, the transfer operator $\phi \mapsto L_a \phi$
as a map from $C^{k+r}(\mathcal{M})$ into $C^k(\mathcal{M})$
is a $C^{r+1}$ map of $a$ and $\phi$.
\end{proposition}

\noindent {\sc Proof.}
By Proposition~\ref{Lcontinuous}, $L_a$ depends continuously on $a$.
For the derivative of $L_a \phi$ with respect to $a$ we find an expression
similar to (\ref{derivative}),
\begin{equation}
\frac{\partial}{\partial a} L_a \phi (x) =
\int_{V_{x,a}} \frac{\partial}{\partial a} k_a(y,x) \phi(y) dm(y)
+\int_{\partial V_{x,a} } s_a(y,x) k_a(y,x) \phi(y) dS(y)
\end{equation}
for some smooth function $s_a$.
It follows that
for $\phi \in C^k(\mathcal{M})$,
$ \frac{\partial}{\partial a} L_a \phi \in C^k(\mathcal{M})$.
This implies differentiability of $(\phi,a) \mapsto L_a \phi$ for $\phi \in C^k(\mathcal{M})$.
Higher differentiability is treated similarly.
\qed\\

The operator $L_a: C^k(\mathcal{M}) \to C^k(\mathcal{M})$ does not depend $C^2$ on $a$, since
$ \frac{\partial^2}{\partial a^2} L_a \phi$
may not exist if $\phi \in C^0(\mathcal{M})$
and $ \frac{\partial^2}{\partial a^2} L_a \phi$ is a $C^{k-1}$ function if $\phi \in C^k(\mathcal{M})$.
What does hold is that $(x,a) \mapsto L_a \phi_a (x)$ is $C^{k+1}$
if $(x,a) \mapsto \phi_a(x)$ is $C^k$.\\


\noindent {\sc Proof of Theorem~\ref{continuation}.}
Let $W$ be the isolating neighborhood for $\mu_{a_0}$.
For $a$ near $a_0$, $f_a(W;\Delta)$ is strictly contained in $W$
and $f_a$ has a unique stationary measure with support in $W$.
Restrict $f_a$ to $W$ for such values of $a_0$.

Consider the transfer operator $L_a$ for $f_a$ acting on
$C^k_0(W)$. Write $F$ for the line in $C^k_0(W)$
spanned by $\phi_{a_0}$. Then $C^k_0(W) = F \oplus
H^k_0(W)$ with $H^k_0(W)$ consisting of $C^k$
functions with vanishing integral;
$$
H^k_0(W) = \{ \phi \in C^k_0(W) \;|\; \int_\mathcal{M} \phi(x) dm(x) =0\}.
$$
Write $\phi_a$ for the eigenvectors of $L_a$ continuing $\phi_{a_0}$ provided
by Proposition~\ref{isolated}.
Decompose $\phi_a = \phi_{a_0} + \psi_a$ with $\psi_a \in H^k_0(W)$.
Then $\psi_a$ is a solution of
$L_a \psi_a  = \psi_a + \phi_{a_0} - L_a\phi_{a_0}$.
Note that $(x,a) \mapsto \phi_{a_0}(x) - L_a\phi_{a_0}(x)$ is $C^\infty$.
The spectrum of $\Res{L_{a_0}}{H^k_0(W)}$ is away from 1.
Proposition~\ref{prop_regular} in Appendix~\ref{sec_ift} implies the result.
\qed

\section{Conditionally stationary measures}\label{sec_escape}

To study average
escape times from open sets we make use of conditionally stationary measures,
which are measures for which on average
a fixed percentage of mass escapes under an iterate.
We recall the notion of conditionally invariant measure,
see \cite{piayor79,pia81,chemartro98,homyou02} for its use
in deterministic dynamics.
Let a map $f : \mathcal{M} \to \mathcal{M}$ be given and restrict $f$ to a domain $W \subset \mathcal{M}$.
Let $V \subset W$ be the set of points in $W$ that are mapped into $W$, points in the
complement of $V$ in $W$ are mapped outside $W$. Consider $f : V \to W$.
A conditionally invariant measure for $f$ on $W$
is a measure $\mu$ on $\mathcal{M}$ so that $\mu(A) = \mu(f^{-1}(A))/ \mu(f^{-1}(W))$
for Borel sets $A \subset W$.

\begin{definition}
Let $f \in R^\infty(\mathcal{M})$. Let $W$ be an open domain in $\mathcal{M}$.
A measure $\bar{\mu}$ on $\overline{W}$ is a
{\em conditionally stationary measure} if
$$
\bar{\mu} (A) = \left. \int_W  P(x,A) d\bar{\mu}(x) \right/ \int_W P(x,W) d\bar{\mu}(x)
$$
for Borel sets $A \subset W$.
\end{definition}

See \cite{ferkesmarpic95,laspea01} where this notion is called a
quasistationary measure.
Note that a conditionally stationary measure is a stationary measure if
$\int_W P(x,W) d\bar{\mu}(x) = 1$, that is, if the support of the
conditionally stationary measure lies inside
$\overline{W}$.

\begin{lemma}
A measure $\bar{\mu}$ on $\overline{W}$ is a
conditionally stationary measure for $f$ if and only if
$\bar{\mu} \times \nu^\infty$ is a conditionally invariant measure
for $(f,\vartheta)$ on $\overline{W} \times \Delta^{\mathbb N}$.
\end{lemma}

\noindent {\sc Proof.} Write $D^1 (x) = \{ {\boldsymbol \omega}
\in \Delta^{\mathbb N} \;|\;
   f(x;{\boldsymbol \omega}) \in W \}$.
Consider $S : \cup_{x \in W}( \{x\} \times
    D^1(x)) \to W \times \Delta^{\mathbb N}$,
    $S(x,{\boldsymbol \omega}) = (f(x; {\omega}) , \vartheta {\boldsymbol \omega})$.
We must show that
the following two statements are equivalent.
\begin{itemize}

\item[$(i)$]
$\bar{\mu}\times \nu^\infty (S^{-1}(A))\left/
 \bar{\mu}\times \nu^\infty (S^{-1}(W \times \Delta^{\mathbb N}))\right. =
            \bar{\mu} \times \nu^{\infty} (A)$
for Borel sets $A \subset W \times \Delta^{\mathbb N}$.

\item[$(ii)$]
$\int_W \int_{\Delta} 1_U (f(x;\omega))
d \bar{\mu}(x) d \nu (\omega)/
\int_W \int_{\Delta} 1_W (f(x; \omega))
d \bar{\mu}(x) d \nu ( \omega)
 = \int_W 1_U(x) d\bar{\mu} (x)$ for Borel sets $U \subset W$.

\end{itemize}
Take a Borel set $U \times V$ with $U \subset W$ and
$V \subset \Delta^{\mathbb N}$ and
compute
\begin{eqnarray}
\nonumber
\bar{\mu}\times \nu^\infty (S^{-1} (U \times V)) & = &
\bar{\mu}\times \nu^\infty  \left( \bigcup_{\omega \in \Delta} f^{-1}(U;\omega) \times \{\omega\}
\times V \right)
\\
\nonumber
& = &  \bar{\mu}\times \nu  \left( \bigcup_{\omega \in \Delta} f^{-1}(U;\omega) \times \{\omega\}
        \right)  \nu^\infty (V)
\\
\label{1}
& = &
\int_W \int_\Delta  1_U (f(x;\omega)) d \bar{\mu}(x) d \nu(\omega) \nu^\infty (V)
\end{eqnarray}
Further
\begin{eqnarray}\label{2}
\bar{\mu} \times \nu^\infty (U \times V) & = & \int_W 1_U (x) d\bar{\mu}(x) \nu^\infty (V)
\end{eqnarray}
Equations (\ref{1}) and (\ref{2})
contain the implication $(i) \Rightarrow (ii)$
when applied for
$U \times \Delta^{\mathbb N}$ for the enumerator and
for $W \times \Delta^{\mathbb N}$ for the denominator.

To show that $(ii)$ implies $(i)$,  note that  (\ref{1}) and (\ref{2}) show that
$(i)$ holds for Borel sets $A$ of the form $U \times V$ if $(ii)$ is assumed.
Therefore it holds for
all Borel sets in $W \times \Delta^{\mathbb N}$.
\qed\\

We continue with the introduction of transfer operators whose fixed points are
the densities of conditionally stationary measures.
The transfer operator $\bar{L}$,
defined for functions in
${\cal L}^1(\overline{W})$ with integral 1, is given by
\begin{equation}
\bar{L} (\phi)  = 1_{\overline{W}} L \phi    \left/ \int_W L \phi(x) dm(x). \right.
\end{equation}
Write
$\tilde{L} \phi = 1_{\overline{W}} L \phi.$


\begin{proposition}\label{prop_barLcompact}
$\tilde{L}$ maps $C^0 (\overline{W})$ into itself and is a compact operator
on it.
\end{proposition}

\noindent {\sc Proof.}
If $k(x,y)$ denotes the density of the stochastic transition function $P(x,\cdot)$, then
$$
\tilde{L} \phi (x) =
  \int_{\overline{W} \cap V_x} k(y,x) \phi(y) dm(y).
$$
Note that $\overline{W}\cap V_x$ depends continuously on $x$.
Recall from
the proof of Proposition~\ref{compact} that we must show that
\begin{itemize}
\item for all $x \in \overline{W}$,
$\{ |\tilde{L}\psi(x)| \;|\; \psi \in B^0(\overline{W}) \}$
is bounded,
\item $\tilde{L}F$ is equicontinuous.
\end{itemize}
Here $B^0(\overline{W})$ is the unit sphere in $C^0(\overline{W})$.
The first item follows as before: $|\tilde{L}\psi(x)| \le \int_{W} k(y,x)dm(y)$  is bounded
by a continuous function and thus bounded.
For the second item we must show that for each $\epsilon >0$ there is $\delta >0$
so that for all $x \in \overline{W}$, $\psi \in B^0(\overline{W})$,
 $|\tilde{L}\psi (x + h) - \tilde{L}\psi (x)| < \epsilon$ if $|h| < \delta$.
Recall, see(\ref{continuityL}),
\begin{eqnarray*}
 \tilde{L}\psi (x + h) - \tilde{L}\psi (x)
& = &
 \int_{W \cap V_{x+h}\cap V_{x}} ( k(y,x+h) - k(y,x) ) \psi (y) dm(y)
\\
& &
+
 \int_{W \cap V_{x+h} \backslash (W \cap V_{x+h}\cap V_{x})} k(y,x+h) \psi(y) dm(y)
\\
& &
-
\int_{W \cap V_{x} \backslash (W \cap V_{x+h}\cap V_{x})} k(y,x) \psi(y) dm(y).
\end{eqnarray*}
Now
$$|\int_{W \cap V_{x+h}\cap V_{x}} ( k(y,x+h) - k(y,x) ) \psi (y) dm(y)| \le
\int_{W \cap V_{x+h}\cap V_{x}} |k(y,x+h) - k(y,x) | dm(y),$$
which is small for $|h|$ small by uniform continuity of $k$.
And
$$
|\int_{W \cap V_{x+h} \backslash (W \cap V_{x+h}\cap V_{x})} k(y,x+h) \psi(y) dm(y)|
\le
\int_{W \cap V_{x+h} \backslash (W \cap V_{x+h}\cap V_{x})} | k(y,x+h)|  dm(y)
$$
is small for $|h|$ small by boundedness of $k$ and uniform continuity
of the volume of $V_x$ in $x$. Similarly for the third term.
This proves equicontinuity.
\qed\\

Recall from Proposition~\ref{Lcontinuous}
that $L$ depends continuously on the random diffeomorphism.
The same argument shows that $\tilde{L}$ depends continuously on
the random diffeomorphism.

\begin{proposition}\label{Ltildecontinuous}
The transfer operator $\tilde{L}$ as a linear map on $C^0(\mathcal{M})$
depends continuously on $f \in R^k(\mathcal{M})$.
\end{proposition}

We obtain conditionally stationary measures by a perturbation argument, perturbing
from an invariant measure.
We do not develop general existence results for conditionally stationary measures,
as such general results are not needed for our purposes.
Let $\{f_a\}$ be a family of random diffeomorphisms depending on a real parameter
$a \in I$. Consider, for $a_0 \in I$, a stationary density  $\phi_{a_0}$ with
support $E_{a_0}$. Let $W$ be a neighborhood of $E_{a_0}$ disjoint from
the supports of possible other stationary densities of $f_{a_0}$.

\begin{proposition}\label{prop_csm}
For $a$ close to $a_0$, $f_{a}$ possesses a
conditionally stationary density $\bar{\phi}_a$ on $W$,
with $\bar{\phi}_{a_0} = \phi_{a_0}$ and
$(x,a) \mapsto \bar{\phi}_a(x)$ continuous in $(x,a)$.
One has
$$
\tilde{L} \bar{\phi}_a = \alpha (a) \bar{\phi}_a,
$$
where $\alpha(a) = \int_W L_a \bar{\phi}_a (x) dm(x)$ depends continuously on $a$.
\end{proposition}

\noindent {\sc Proof.}
The operator $\tilde{L}_a$ varies continuously with $a$ and therefore
possesses a single eigenvalue
close to 1 for
$a$ close to $a_0$.
The function $\bar{\phi}_a$ is the corresponding
eigenfunction.
\qed\\

Recall the definition of the escape time $\chi_a (x,{\boldsymbol
\omega})$ for $x \in W$ and ${\boldsymbol \omega} \in
\Delta^{\mathbb N}$:
\begin{equation*}
\chi_a(x,{\boldsymbol \omega}) = \min \{k \;|\; f^k_a
(x;{\boldsymbol \omega}) \not \in W\}.
\end{equation*}

\begin{lemma}\label{lemma_averageexit}
$$
\int_W \int_{\Delta^{\mathbb N}} \chi_a (x ;{\boldsymbol \omega})
d\nu^\infty(\omega) \bar{\phi}_a (x) dm(x)
 = \frac{1}{1-\alpha(a)}.
$$
\end{lemma}

\noindent {\sc Proof.}
%
Let $S_a(x,{\boldsymbol \omega}) = (f_a(x; {\omega}) , \vartheta
{\boldsymbol \omega})$. Write
$$ D^n_a = \{ (x,{\boldsymbol \omega}) \;|\; f^n_a(x;\omega_1,\ldots,\omega_n) \in W \}$$
 for the set of points in $W \times \Delta^{\mathbb N}$ that remain in
 $W \times \Delta^{\mathbb N}$ for $n$ iterates of $S_a$.
 The exit set $E^n_a$ of points that leave $W \times \Delta^{\mathbb N}$
 in $n$ iterates equals $D^{n-1}_a \backslash D^{n}_a$.
 Thus $\chi_a (x,{\boldsymbol \omega}) = n$ on $E^n_a$.
Write $\bar{\mu}_a$ for the conditionally stationary measure with density $\bar{\phi}_a$. From (\ref{1})
with $A = W \times \Delta^{\mathbb N}$ we get
 \begin{eqnarray*}
\alpha(a) & = & \int_W P(x,W) d \bar{\mu}_a(x)
\\
& = & \int_W \int_{\Delta} 1_W(f(x;\omega)) d \nu(\omega) d\bar{\mu}_a (x)
\\
& = & \bar{\mu}_a\times \nu^\infty (S_a^{-1} (W \times \Delta^{\mathbb N})).
\end{eqnarray*}
It follows that
$
\bar{\mu}_a \times \nu^\infty (E^k_a) = \alpha^{k-1} - \alpha^k = \alpha^{k-1}(1-\alpha).
$
Calculate
\begin{eqnarray*}
\int_W \int_{\Delta^{\mathbb N}} \chi_a (x ;{\boldsymbol \omega})
d\nu^\infty(\omega) d \bar{\mu}_a (x)
& = &
\sum_{k=1}^\infty \bar{\mu}_a \times \nu^\infty (E^n_a)
\\
& = & \sum_{k=1}^\infty k \alpha^{k-1} (1-\alpha)
\\
& = & \frac{1}{1-\alpha}.
\end{eqnarray*}
\qed\\

As a corollary we obtain that
the average escape time from $W$ goes to infinity as $a \to a_0$.
More precise estimates are derived in the following section.

\begin{proposition}\label{theorem_exittimes_continuous}
 $\int_{W} \int_{\Delta^{\mathbb N}} \chi_a (x ,{\boldsymbol \omega}) d\nu^\infty({\boldsymbol \omega}) dm(x) $
converges to $\infty$ as $a \to a_0$.
\end{proposition}

\noindent {\sc Proof.}
Let $E_a$ be the interior of the support of $\bar{\phi}_a$.
\begin{eqnarray*}
\int_{W} \int_{\Delta^{\mathbb N}} \chi_a (x ,{\boldsymbol
\omega})
              d\nu^\infty({\boldsymbol \omega}) dm(x)
& \ge & \int_{E_a} \int_{\Delta^{\mathbb N}} \frac{\chi_a (x
,{\boldsymbol \omega})}{\bar{\phi}_a(x)}  d\nu^\infty({\boldsymbol
\omega}) \bar{\phi}_a(x) dm(x)
\\
& \ge & C \int_{E_a} \int_{\Delta^{\mathbb N}} \chi_a (x
,{\boldsymbol \omega}) d\nu^\infty({\boldsymbol \omega})
\bar{\phi}_a(x) dm(x)
\\
& = & \frac{C}{1-\alpha(a)},
\end{eqnarray*}
for $C = 1/\max \{x \in W \;|\; \bar{\phi}_a (x) \}$.
For $x$ from the support $E_{a_0}$, the image $f_{a_0}(x;\Delta)$ is contained
in $E_{a_0}$ by invariance of $E_{a_0}$.
By continuity of $f_a$,  $f_{a}(x;\Delta) \subset W$ for $a$ close enough to $a_0$.
Since $\bar{\phi}_a$ depends continuously on $a$,
 $1 - \alpha(a) \le \int_{W\backslash E_{a_0}} \bar{\phi}_a(x)dm(x)$
converges to $0$ as $a \to a_0$.
The proposition follows.
\qed

\section{Escape times}\label{sec_escapeII}

In this section estimates for the average escape time from small neighborhoods of the support
of a stationary measure that undergoes a bifurcation are derived,
as function of the unfolding parameter.

Let $\{ (f_a, g_a)\}$, $a \in I$, be a smooth
one parameter family of random diffeomorphisms on $\mathcal{M}$.
Suppose that $a_0 \in I$ is a bifurcation value for $\{ (f_a, g_a)\}$.
Let $\mu$ be a stationary measure for $f_{a_0}$ involoved in a 
bifurcation.
Let $W$ be a small neighborhood of $E$.
By Proposition~\ref{prop_csm},
for $W$ sufficiently close to $E$ and $a$ near $a_0$,
$\{ (f_a, g_a)\}$ possesses a unique conditionally stationary measure $\bar{\mu}_a$ with support in $W$.
Write $\bar{\phi}_a$ for the density of $\bar{\mu}_a$.
The transfer operator $\bar{L}_a$ acting on $C^0(\overline{W})$, 
has $\bar{\phi}_a$ as a unique fixed point.

Let $X^0$ be the set of points in $\overline{W}$ with
$\partial V_{x,a} \cap \partial W \neq \emptyset$ for $x \in X^0$.
For $i \ge 0$, define $X^{i+1} = f(X^i ; \partial \Delta)$.
We suppress the dependence of $X^i$ on $a$ from the notation.

\begin{lemma}
If $x_i \in X^i$, then for $x_{i+1} \in f(x_i;\partial \Delta)$
one has $x_i \in \partial V_{x_{i+1},a}$.
\end{lemma}

\noindent {\sc Proof.} This is clear from the definition. \qed \\


At $x \in X^0$, the boundary of
$V_{x,a} \cap W$ varies continuously  but not smoothly with $(x,a)$.
It follows that $\bar{L}_a (\phi)$ cannot be expected to be
more than continuous on $X^0$ even for smooth $\phi$.

\begin{lemma}\label{lem_smoothoutside}
Suppose that $(x,a)  \mapsto \phi_a(x)$ is $C^k$ outside $X^0 \cup \cdots \cup X^{k-1}$,
such that derivatives up to order $k$ are bounded and their restrictions
to a component of $W \backslash (X^0 \cup \cdots \cup X^{k-1})$
extend continuously to
the boundary of the component.
Then $(x,a) \mapsto \bar{L}_a \phi_a(x)$ is $C^{k+1}$
outside $X^0 \cup \cdots \cup X^{k}$.
Likewise, derivatives up to order $k+1$ are bounded and their restrictions
to a component of $W \backslash (X^0 \cup \cdots \cup X^{k})$
extend continuously to the boundary of the component.
\end{lemma}

\noindent {\sc Proof.}
For $x \not \in X^0$, the derivative of $\tilde{L}_a \phi$ is of the form
\begin{equation}\label{derivative_capW}
D (\tilde{L}_a \phi) (x) = \int_{V_{x,a} \cap W} \frac{\partial}{\partial x} k_a(y,x) \phi(y) dm(y)
+ \int_{\partial (V_{x,a} \cap W)} n_a(y,x) k_a(y,x) \phi(y) dS(y).
\end{equation}
for a smooth function $k_a $ and a piecewise smooth function $n_a$
(smooth outside
the intersection of $\partial W$ with $\partial V_{x,a}$).
This identity shows that $\tilde{L}_a \phi$ is $C^1$ with
bounded derivatives outside $X^0$,
for any $\phi \in C^0(\overline{W})$. The same holds for $\bar{L}_a (\phi)$.
Higher order derivatives
are treated inductively.
Similar to (\ref{derivative_capW})
 one has
\begin{equation}\label{derivative_acapW}
\frac{\partial}{\partial a} \tilde{L}_a \phi (x) =
\int_{V_{x,a} \cap W} \frac{\partial}{\partial a} k_a(y,x) \phi(y) dm(y)
+\int_{\partial (V_{x,a} \cap W)} s_a(y,x) k_a(y,x) \phi(y) dS(y)
\end{equation}
for
some piecewise smooth function $s_a$.
This shows that $(x,a)\mapsto \tilde{L}_a \phi_a(x)$ is $C^1$ outside $X^0$
for continuous functions $(x,a)\mapsto \phi_a(x)$.
The derivatives of $(x,a)\mapsto \tilde{L}_a \phi_a(x)$
on $W \backslash X^0$ are bounded; moreover
the derivatives on a component of $W \backslash X^0$ extend 
continuously to the boundary of the component.

Higher order derivatives are treated inductively.
Suppose that $(x,a)  \mapsto \phi_a(x)$ is $C^k$ outside $X^0 \cup \cdots \cup X^{k-1}$,
such that derivatives up to order $k$ are bounded and their restrictions
to a component of $W \backslash (X^0 \cup \cdots \cup X^{k-1})$
extend continuously to
the boundary of the component.
Then $(x,a) \mapsto \tilde{L}_a \phi_a(x)$ is $C^{k+1}$
outside $X^0 \cup \cdots \cup X^{k}$.
Likewise, derivatives up to order $k+1$ are bounded and their restrictions
to a component of $W \backslash (X^0 \cup \cdots \cup X^{k})$
extend continuously to
the boundary of the component.

The transfer operator $\bar{L}_a$ is the composition of the linear map
$\tilde{L}_a$ and the projection
$$
\Pi ( \phi)  =  \phi/\int_W \phi(x)dm(x).
$$
The projection $\Pi$ is a smooth map which is well defined near $\phi_{a_0}$
in $C^0(\overline{W})$, a direct computation shows
$$
D\Pi (\phi) h = \frac{1}{\int_W \phi(x)dm(x)} h - \frac{\phi}{(\int_W \phi(x)dm(x))^2} \int_W h(x)dm(x).
$$
Also,
$$\frac{\partial^i}{\partial a^i} \int_W \phi_a (x) dm(x) =
\int_W \frac{\partial^i}{\partial a^i} \phi_a (x) dm(x).
$$
This implies that also $(x,a) \mapsto \bar{L}_a (\phi_a) (x)$ is $C^{k+1}$
outside $X^0 \cup \cdots \cup X^{k}$ and has bounded derivatives.
\qed

\begin{proposition}\label{prop_barphi}
For each $k \ge 1$, $\bar{\phi}_a$ is
$C^k$ outside $X^0 \cup \cdots \cup X^{k-1}$ jointly in $(x,a)$,
the derivatives up to order $k$
are uniformly bounded.
\end{proposition}

\noindent {\sc Proof.}
Proposition~\ref{prop_csm} gives that $\bar{\phi}_a(x)$ is continuous in $(x,a)$.
Recall that the support $E$ of $\mu$ consists of finitely many, say $k$,
connected components, permuted cyclically
by the random diffeomorphism. An iterate of $f_{a_0}$ thus maps
each component into itself.
The restriction of the transfer operator $L^k_{a_0}$ to a small neighborhood of
$E$ has a single eigenvalue 1 and
a remaining spectrum strictly inside the unit circle \cite{gay01}.
There is therefore no loss in generality to assume that the support of
$\mu$ consists of a single connected component $E$, which we will assume
for the remainder of the proof.

Write $H^k(\overline{W}) = \{\psi \in C^k(\overline{W}) \;|\; \int_W \psi = 0\}$.
Define the operator
$T_a : H^0(\overline{W}) \to H^0 ( \overline{W})$
by
\begin{equation}
T_a (\psi) = \bar{L}_a ( \phi_{a_0} + \psi) - \phi_{a_0}.
\end{equation}
Decompose $\bar{\phi}_a = \bar{\phi}_{a_0} + \bar{\psi}_a$,
so that $T_a (\bar{\psi}_a) = \bar{\psi}_a$. From the proof of Lemma~\ref{lem_smoothoutside}, we get that 
$T_a$ a smooth map on $H^0(\overline{W})$;
$$
DT_a (\psi) = D\Pi (\tilde{L}_a (\phi_{a_0} + \psi)) \tilde{L}_a.
$$
However, $\tilde{L}_a$ maps continuously differentiable functions to
continuous functions, so that $T_a$ does not define a map from
$H^k( \overline{W})$, $k\ge 1$, to itself.
As a consequence we cannot obtain smoothness properties of $\bar{\phi}_a$
by applying the implicit function theorem.
To get smooth dependence of $\bar{\phi}_a$ outside sets $X^i$
we reason as follows.
We derive equations the derivatives of $\bar{\phi}_a$ must
satisfy, establish that the equations
can be solved, and show that the solutions are the derivatives
of $\bar{\phi}_a$.
The reasoning follows the lines of the proof of Proposition~\ref{prop_regular}
in Appendix~\ref{sec_ift}.


To prove that $\bar{\psi}_a$ varies $C^1$ with $a$ in points outside $X^0$,
note that $\frac{\partial}{\partial a} \bar{\psi}_a (x)$ should be a solution $M_a(x)$
 to
\begin{equation}\label{Ma}
 \frac{\partial}{\partial a} T_a (\bar{\psi}_a) (x) + D  T_a (\bar{\psi}_a)  M_a(x)
 = M_a(x)
\end{equation}
We claim that this equation is uniquely solvable.
The spectral radius of $DT_{a_0}(0) = L_{a_0}$ is smaller than 1.
As a consequence of the continuous dependence of $\tilde{L}_a$ on $a$
(see Proposition~\ref{Ltildecontinuous}), also
$D T_a (\psi_a)$ varies continuously with $a$.
For $a$ sufficiently close to $a_0$, the
spectral radius of $D T_a  (\bar{\psi}_a)$ is therefore also
smaller than 1.
Hence
\begin{equation}\label{series}
\left( I - D T_a (\bar{\psi}_a) \right)^{-1} =   I + \sum_{i=1}^\infty (D T_a (\bar{\psi}_a))^i,
\end{equation}
see \cite{kat66}.
This formula can be applied for $D T_a (\bar{\psi}_a)$ acting on
${\cal L}^2$ functions.
Indeed,
$DT_a$ is compact on
${\cal L}^2 (\overline W)  \subset {\cal L}^1 (\overline W)$,
see Remark~\ref{remark_l2}, and has spectrum strictly inside
the unit circle in ${\mathbb C}$. From
\begin{eqnarray}
M_a(x) & = & \nonumber
\left( I - D T_a (\bar{\psi}_a) \right)^{-1}
\frac{\partial}{\partial a} T_a  (\bar{\psi}_a) (x)
\\
\label{whatitisreturn}
& = &
\left(  I + D T_a (\bar{\psi}_a) + (DT_a (\bar{\psi}_a))^2 + (D T_a (\bar{\psi}_a))^3 + \cdots \right)
\frac{\partial}{\partial a} T_a  (\bar{\psi}_a) (x),
\end{eqnarray}
we get that $M_a$ is continuous outside $X^0$ since it
equals the sum of $\frac{\partial}{\partial a} T_a  (\bar{\psi}_a)$
and a uniform limit of continuous functions (compare Lemma~\ref{lem_smoothoutside}).
In particular $M_a$ is uniformly bounded and has continuous extensions to the closure
of components of $W \backslash X^0$.
We must show that
$$
\left| \bar{\psi}_{a+h} (x) - \bar{\psi}_a (x)  - M_a (x)h \right| = o(|h|)
$$
for $x \not \in X^0$, as $h \to 0$.
Consider $\gamma_a (x) = \bar{\psi}_{a+h} (x) - \bar{\psi}_a (x)$ for $x \not \in X^0$.
Then
\begin{eqnarray}\nonumber
\gamma_a (x) & = & T_{a+h} (\bar{\psi}_a + \gamma_a)(x) - T_a (\bar{\psi}_a) (x)
\\
\label{psia}
& = &
D T_a (\bar{\psi}_a)\gamma_a (x) + \frac{\partial}{\partial a} T_a ( \bar{\psi}_a) (x) h + R(x),
\end{eqnarray}
where
$$
R(x)  =  T_{a+h} (\bar{\psi}_a + \gamma_a)(x) - T_a (\bar{\psi}_a) (x)
-  D T_a (\bar{\psi}_a) \gamma_a (x) - \frac{\partial}{\partial a} T_a (\bar{\psi}_a) (x) h.
$$
We claim that
for any $\epsilon >0$ there is $\delta>0$ so that
$|R| < \epsilon ( |\gamma_a | + |h|)$,  if
$|h|$ and $|\gamma_a|$ are smaller than $\delta$.
Since $\gamma_a (x)$ is continuous in $h$ we may
further restrict $\delta$ in this estimate so that
$|R| < \epsilon ( |\gamma_a| + |h|)$ holds for
$|h|$ smaller than $\delta$.
Further, $\left( I - D T_a (\bar{\psi}_a)  \right) \gamma_a (x)
   =  \frac{\partial}{\partial a} T_a( \bar{\psi}_a) (x) h + R(x)$.
Using (\ref{series}) and the bound on $|R|$
gives
$|\gamma_a|  \le k |h|$ for some $k$ if $|h| < \delta$.
Therefore $|R| < \epsilon (1+k) |h|$ for some $k>0$, if $|h| < \delta$.
Now (\ref{Ma}) and (\ref{psia}) give
$$
\gamma_a (x) - M_a (x) h  = \left( I - D  T_a (\bar{\psi}_a) \right)^{-1} R (x).
$$
Using (\ref{series}) it follows that
$\left| \gamma_a  - M_a h \right| = o(|h|)$, $h \to 0$.
This proves that $M_a$ equals the
partial derivative $\frac{\partial}{\partial a} \bar{\psi}_a$.



Higher orders of differentiability are proved by induction.
Assume that $(x,a) \mapsto \bar{\psi}_a(x)$
has been shown to be $C^k$ outside $X^0\cup\cdots\cup X^{k-1}$.
Recall from Lemma~\ref{lem_smoothoutside} that
for $C^k$ maps $(x,a) \mapsto \bar{\psi}_a(x)$, $\frac{\partial}{\partial a} \bar{L}_a (\bar{\psi}_a)$
is $C^{k}$ outside  $X^0\cup\cdots\cup X^{k}$.
The right hand side of (\ref{whatitisreturn})
is therefore $C^{k}$ outside  $X^0\cup\cdots\cup X^{k}$.
The above reasoning shows that $M_a = \frac{\partial}{\partial a} \bar{\psi}_a$
outside  $X^0\cup\cdots\cup X^{k}$.
Therefore
$\frac{\partial}{\partial a} \bar{\psi}_a$
is $C^k$  outside  $X^0\cup\cdots\cup X^{k}$.
Also $D\bar{\psi}_a = D  (T_a (\bar{\psi}_a))$ is $C^k$ outside  $X^0\cup\cdots\cup X^{k}$,
so that
$(x,a) \mapsto \bar{\psi}_a(x)$ is $C^{k+1}$ outside $X^0\cup\cdots\cup X^{k}$.
The same clearly holds for $(x,a) \mapsto \bar{\phi}_a(x)$.
\qed\\

\noindent {\sc Proof of Theorem~\ref{theorem_exittimes_smooth}.}
We repeat the computation in the proof of Proposition~\ref{theorem_exittimes_continuous}.
Let $E_a$ be the interior of the support of $\bar{\phi}_a$.
Applying Lemma~\ref{lemma_averageexit},
\begin{eqnarray*}
\int_{W} \int_{\Delta^{\mathbb N}} \chi_a (x ,{\boldsymbol
\omega})
              d\nu^\infty({\boldsymbol \omega}) dm(x)
& \ge & \int_{E_a} \int_{\Delta^{\mathbb N}} \frac{\chi_a (x
,{\boldsymbol \omega})}{\bar{\phi}_a(x)}  d\nu^\infty({\boldsymbol
\omega}) \bar{\phi}_a(x) dm(x)
\\
& \ge & C \int_{E_a} \int_{\Delta^{\mathbb N}} \chi_a (x
,{\boldsymbol \omega}) d\nu^\infty({\boldsymbol \omega})
\bar{\phi}_a(x) dm(x)
\\
& = & \frac{C}{1-\alpha(a)},
\end{eqnarray*}
for $C = 1/\max \{x \in \overline{W} \;|\; \bar{\phi}_a (x) \}$.
By Proposition~\ref{prop_barphi}, $(x,a) \mapsto \bar{\phi}_a(x)$
is $C^k$ almost everywhere
and has uniformly bounded derivatives.
For each integer $k$ there is a constant $C$ with
$|\bar{\phi}_a| \le C |a-a_0|^k$ on $\overline{W \backslash E_{a_0}}$.
As in the proof of Proposition~\ref{theorem_exittimes_continuous} we
get that for each $k$ there is a constant $C_k>0$,
so that $1 - \alpha (a)  \le C_k |a-a_0|^k$.
\qed

\section{Decay of correlations}\label{sec_correlations}

Consider a random family $\{f_a\}$ restricted to an isolating neigborhood $W$
of a stationary measure $\mu_a$, for all values of $a$ form an interval $I$.
The transfer operator $L_a$ on $C^k_0(W)$ possesses a single eigenvalue at 1.
If the support of $\mu_a$ consists of $r$ components,  $L_a$ has eigenvalues
$e^{2\pi i/j}$, $0 \le j < r$, on the unit circle in the complex plane.
These eigenvalues make up the peripheral spectrum of $L_a$, see Remark~\ref{remark_tobias}.
In this section we consider bifurcations in which the number
of components of the support of $\mu_a$ changes.
We will see how the rate of decay of correlations varies with the parameter $a$,
providing a proof of Theorem~\ref{theorem_corr}.


\begin{proposition}\label{lambdasmooth}
Let $\{f_a\}$, $a\in I$,  
be a family of random diffeomorphisms with an isolating neighborhood $W$.
The eigenvalues and eigenvectors of the peripheral spectrum of $L_a$ on $C^k_0(W)$ 
vary smoothly with $a$.
\end{proposition}

\noindent {\sc Proof.} Let $\lambda_a$ be an eigenvalue that
depends continuously on $a$ and lies on the unit circle for
$a=a_0$. Since $\mu_{a_0}$ is an isolated stationary measure,
$\lambda_{a_0}$ is a simple eigenvalue 
(see Remark~\ref{remark_tobias}).
Proposition~\ref{prop_regular}
in Appendix~\ref{sec_ift} implies the result. \qed\\

Recall (\ref{P-F}) and Lemma~\ref{lemma_transfer}.
Write
\begin{eqnarray*}
L^n_a \varphi (x) & = & \int_{\Delta^n} P_{f_a(x;\omega_1,\ldots,\omega_n)} \varphi (x)
                             d\nu(\omega_1)\cdots d\nu(\omega_n),
\\
U^n_a \psi (x) & = & \int_{\Delta^n} \psi \circ f_a (x;\omega_1,\ldots,\omega_n)
                        d\nu(\omega_1)\cdots d\nu(\omega_n).
\end{eqnarray*}
As in the computation for  Lemma~\ref{lemma_transfer},
\begin{equation}\label{adjoint}
\int_\mathcal{M} L^n_a \varphi (x) \psi(x) dm(x) = \int_\mathcal{M} \varphi(x) U^n_a \psi(x) dm(x).
\end{equation}
After these preparations we now prove the statements
on the speed of decay of correlations. First consider a single random map $f$.\\

\noindent {\sc Proof of Proposition~\ref{prop_corr}.}
Let $\phi$ be the stationary density.
Write $L^n \varphi =  \left( \int_\mathcal{M} \varphi(y) dm(y) \right) \phi + R^n \varphi$.
Compute
\begin{eqnarray*}
\int_\mathcal{M} \varphi(x) U^n \psi(x) dm(x) &=& \int_\mathcal{M} L^n \varphi(x) \psi(x) dm(x)
\\
&=& \int_\mathcal{M} \left[ \left( \int_\mathcal{M} \varphi(y) dm(y) \right) \phi (x) + R^n \varphi(x) \right] \psi(x) dm(x),
\end{eqnarray*}
so that
\begin{eqnarray*}
\left| \int_\mathcal{M} \varphi(x) U^n \psi(x) dm(x) - \int_\mathcal{M} \varphi(x) dm(x) \int_\mathcal{M} \psi(x) \phi(x) dm(x) \right|
&=&
\left| \int_\mathcal{M} R^n\varphi(x) \psi(x) dm(x)  \right|.
\end{eqnarray*}
Note that the spectral radius of $R$ is smaller than 1. 
By continuity of $R$, there is $N>0$ so that $\| R^N \| < 1$
for all $a$ near $a_0$. 
Hence for $n \in {\mathbb N}$, 
$\|  R^n \| < C \eta^n$
for some $C>0, \eta < 1$.
The proposition follows from
\[
\left| \int_\mathcal{M} R^n\varphi(x) \psi(x) dm(x)  \right|
\le \| R^n \varphi\|_{{\cal L}^2(\mathcal{M})} \| \psi\|_{{\cal L}^2(\mathcal{M})}
\le C \eta^n  \|\varphi\|_{{\cal L}^2(\mathcal{M})} \| \psi\|_{{\cal L}^2(\mathcal{M})}.
\]
\qed\\

\noindent {\sc Proof of Theorem~\ref{theorem_corr}.}
This is proved by following the computation in the proof
of Proposition~\ref{prop_corr} above 
and noting that $L_a$ has for $a>a_0$
a single eigenvalue $1$ and $k-1$ eigenvalues that have moved smoothly into
the unit circle.
Write $\eta_a$ for the largest radius of the eigenvalues of $L_a$ that lie
inside the unit circle. As a consequence of the smooth dependence of the eigenvalues
near the unit circle, see Proposition~\ref{lambdasmooth},
$\eta_a$ is a smooth function of $a$.

We claim that there exists $C>0$ so that for all $a$ near $a_0$,  $\|R^n\| \le C \eta_a^n$.
For $a=a_0$, let $E$ be the union of the
eigenspaces for the eigenvalues in the peripheral spectrum.
For $a$ near $a_0$, let $E_a$ be the continuation of $E_{a_0} = E$.
As $E_a$ is finite dimensional and has a basis depending smoothly on $a$,
it is clear that there exists $C>0$ so that for $a>a_0$, $n \in {\mathbb N}$,
\begin{equation}\label{res}
\| \Res{R^n}{E_a} \| \le C \eta_a^n.
\end{equation}
Let $F$ be a subspace of ${\cal L}^2(W)$ complementary to $E$.
Write $P_a$ for the projection to $E_a$ along $F$.
Then $R = P_a R + (I-P_a) R$. 
By continuity of $R$ and $P_a$, there is $N>0$ so that $\| \left( (I-P_a) R \right)^N \| <1$
for all $a$ near $a_0$. 
Hence for $n \in {\mathbb N}$,
\begin{equation}\label{notres}
\| \left( (I-P_a) R \right)^n \| < C \nu^n
\end{equation}
for some $C>0, \nu < 1$.
Now (\ref{res}) and (\ref{notres}) prove the claim.
As before, the proposition follows from
\[
\left| \int_\mathcal{M} R^n\varphi(x) \psi(x) dm(x)  \right|
\le \| R^n \varphi\|_{{\cal L}^2(\mathcal{M})} \| \psi\|_{{\cal L}^2(\mathcal{M})}
\le C \eta_a^n  \|\varphi\|_{{\cal L}^2(\mathcal{M})} \| \psi\|_{{\cal L}^2(\mathcal{M})}.
\]
\qed

\section{One dimensional random maps}\label{sec_1D}

The most complete description of bifurcations in
smooth random maps is derived for random maps
in one dimension.
Consider a random endomorphism $f (x;\omega)$
on the circle ${\mathbb S}^1$. The random parameter $\omega$
is drawn from $\Delta = [-1,1]$.
What is proved below for random endomorphisms on the circle
holds with obvious modifications for
random endomorphisms on a compact interval that is mapped
inside itself by all endomorphisms.

A pathological example occurs if $f(x;\omega)$ is constant in $x$;
the (unique) stationary measure is then a push forward of the
measure on $\Delta$. To avoid pathologies we assume
the open and dense condition
that the critical points of each map $x \mapsto f(x;\omega)$
have finite order.
Also under this condition one finds that
the regularity of stationary measures for random endomorphisms
is substantially less then for random diffeomorphisms;
their densities are only continuous.

\begin{theorem}\label{stationary_measures_1D}
The random endomorphism $f \in R^\infty ({\mathbb S}^1)$ possesses
a finite number of ergodic stationary measures
$\mu_1,\ldots, \mu_m$
with mutually disjoint supports $E_i,\ldots, E_m$.
All stationary measures are linear combinations of
$\mu_1,\ldots,\mu_m$.

The support $E_i$ of $\mu_i$ consists of the closure of a finite
number of connected open sets
$C_i^1, \ldots, C_i^p$
that are moved cyclically by $f(\cdot ;\Delta)$.
The density $\phi_i$ of $\mu_i$ is a $C^0$ function on ${\mathbb S}^1$.
\end{theorem}

\noindent {\sc Proof.}
The condition on the critical points of $x\mapsto f(x;\omega)$
implies that $V_x$ varies continuously with $x$.
The reasoning used to prove Proposition~\ref{compact} shows that
the transfer operator $L$ maps ${\cal L}^1({\mathbb S}^1)$ into $C^0({\mathbb S}^1)$.
This implies continuity of invariant densities.
\qed

\begin{theorem}\label{isolated_1D}
Let $\mu$ be an isolated ergodic stationary measure of $f \in R^\infty ({\mathbb S}^1)$
with density $\phi$.
Then each $\tilde{f} \in R^\infty ({\mathbb S}^1)$ sufficiently close to
$f$ possesses a unique ergodic stationary measure $\tilde{\mu}$
with support
in $V$. The density $\tilde{\phi}$ of $\tilde{\mu}$ is $C^0$ close
to $\phi$.
\end{theorem}

\noindent {\sc Proof.} As in the proof of Theorem~\ref{isolated}.
Note that $L$ is a compact operator on $C^0(\mathcal{M})$, compare the proofs
of Theorems~\ref{compact} and \ref{prop_barLcompact}.\qed\\

Recall that iterates of a random map $f$
are defined through (\ref{iterate}).
A periodic point $\bar{x}$ of period $k$ is a point satisfying
$f^k(\bar{x};\omega_1,\ldots,\omega_k) = \bar{x}$ for some
$\omega_1,\ldots,\omega_k \in \Delta$.
It is hyperbolic if
$\frac{d}{dx} f^k(x;\omega_1,\ldots,\omega_k)$ at $x=\bar{x}$
differs from $0, 1, -1$. By the implicit function theorem,
a family $\{f_a\}$ of random endomorphisms with
$f_{a_0} = f$ possesses
a hyperbolic periodic point
$\bar{x}_a$, $\bar{x}_{a_0} = \bar{x}$,
for $a$ near $a_0$
and for the same values
of $\omega_1,\ldots,\omega_k$,
depending smoothly on $a$.

\begin{theorem}\label{stable_1D}
The set of stable random endomorphisms in $R^\infty ({\mathbb S}^1)$
is open and dense.
\end{theorem}

\noindent {\sc Proof.}
Take $f \in R^\infty ({\mathbb S}^1)$.
If the entire circle is the support of a stationary measure of $f$,
then $f$ is stable by Theorem~\ref{isolated_1D}.
Suppose that $\mu$ is a stationary measure whose support $E$ is a union
$\cup_{i=0}^{k-1} C_i$ of intervals $C_i$ mapped cyclically by $f(\cdot;\Delta)$:
$f(C_i;\Delta) = C_{i+1}$ (the indexes are taken modulo $k$).
If $\mu$ is an isolated measure, $f$ restricted to an isolating neighborhood
of $E$ is stable.

The measure $\mu$ is certainly isolated if for each boundary point $x \in E$,
either
$$
f^k(x;\Delta^{\mathbb N}) \subset \mbox{interior } E,
$$
or
$$
f^k(x;\omega_1,\ldots,\omega_k) = x, \;\;
f^j(x;\omega_1,\ldots,\omega_j)  \in \mbox{interior }  E
$$
for some $\omega_1,\ldots,\omega_k \in \Delta$, $j < k$.
Indeed, invariance of $E$ shows that
in both cases $f^k(y; \Delta^{\mathbb N}) \in E$ for any $y$ near $x$.

If not all boundary points are as above, then there is
a boundary point $x \in \partial E$ so that
$f^l(x;\omega_1,\ldots,\omega_l) = x$ for $l=k$ or $l=2k$
($l$ minimal)  and  $f^j(x;\omega_1,\ldots,\omega_j) \in \partial E$ for $0<j<l$.
Write $x_0=x$ and $x_j = f^j(x;\omega_1,\ldots,\omega_j)$ for $j>0$. From 
$x_j \in \partial E$, $x_{j+1} = f(x_j;\omega_{j+1}) \in \partial E$ and
$\frac{\partial}{\partial \omega} f(\cdot;\omega) \neq 0$,
we see that $\omega_{j+1} \in \partial \Delta$.
Thus $\omega_1,\ldots,\omega_l$ are all contained in $\partial \Delta$.
Note that $\frac{d}{dx} f^l(x;\omega_1,\ldots,\omega_l) \ge 0$ since otherwise
$x$ is an interior point of $E$.

For $f \in R^\infty({\mathbb S}^1)$, there are a neighborhood $U$ of $f$ and
an integer $N$ so that for 
each $\tilde{f} \in U$, the support of the union of its 
stationary measures has
at most $N$ connected components.
A random periodic orbit in the boundary of the support of
a stationary measure of $\tilde{f} \in U$ therefore has its period
bounded by $2 N$.

By transversality techniques  a number of arbitrary small perturbations of 
$f$ are carried through. 
The perturbations affect $f(\cdot;\omega)$ for $\omega \in \partial \Delta$ 
and can be extended to other values of $\omega$ using test functions.
We will not present the detailed perturbations, but
refer to \cite[Section~III.2]{melstr93} for a description of the techniques.
By a small perturbation of $f$ we may assume that the graph of each 
map $f^i (\cdot; (\partial \Delta)^i)$, $1 < i \le 2N$,  
intersects the diagonal in 
${\mathbb S}^1\times {\mathbb S}^1$ transversally. That is, 
\begin{hyp}\label{h_1}the periodic orbits
of period $i \le 2 N$
for $f(\cdot;\partial \Delta)$ are hyperbolic.
\end{hyp}
There is then a bounded number of random periodic orbits with period
bounded by $2 N$.
A further small perturbation ensures that 
\begin{hyp}\label{h_2}each periodic point $x$ of period
$i \le 2 N$ is periodic for only one sequence 
$\omega_1,\ldots,\omega_i \in \partial \Delta$.
\end{hyp}
Write  ${\cal P}$ for the points in these periodic orbits.
Recall that the number of critical points of $f(\cdot;\partial \Delta)$
is finite.
A final small perturbation ensures that 
\begin{hyp}\label{h_3}the critical values of
$f(\cdot;\partial \Delta)$ are disjunct from ${\cal P}$.
\end{hyp}
Conditions~\ref{h_1},~\ref{h_2},~\ref{h_3} are clearly open and thus 
describe an open and dense subset of $R^\infty ({\mathbb S}^1)$.

Consider $f$ from this open and dense set.
Let $\mu$ be a stationary measure of $f$ with support $E$.
Let $x$ be a boundary point of $E$ belonging to a 
periodic orbit in $\partial E$.
By \ref{h_1}, $x$ belongs to a hyperbolic periodic orbit.
By \ref{h_2}, there is a unique graph $f^l(\cdot,\omega_1,\ldots,\omega_l)$
with $\omega_1,\ldots,\omega_l \in \partial \Delta$
through $x$.
It is not possible that
$\frac{d}{dx} f^{l} (x ; \omega_1,\ldots,\omega_l) > 1$, since other orbits
would then be repelled and $x$ would not be in the boundary of $E$.
Hence $0<  \frac{d}{dx} f^{l} (x ; \omega_1,\ldots,\omega_l)<1$: 
the random periodic orbit through $x$ is an attracting periodic orbit
for $f^{l} (\cdot ; \omega_1,\ldots,\omega_l)$.
By \ref{h_3}, there are no interior points in $E$ being mapped onto $x$
under iterates of $f$.
As a consequence,  $\mu$ is isolated. Therefore $f$ is stable.
\qed\\

As a next step we consider one parameter families of random maps
and show that bifurcations typically occur at isolated
parameter families. The theorem below moreover describes
the possible codimension one bifurcations.
The space of smooth families of smooth random maps $x \mapsto f_a (x;\omega)$,
$a \in I$,
will be given the uniform $C^k$ topology as maps
$(x,\omega,a) \mapsto f_a(x;\omega)$ on ${\mathbb S}^1\times \Delta\times I$.

We start with a description of three types of bifurcations
caused by violation of one of the conditions \ref{h_1}, \ref{h_2}, \ref{h_3}.
These are proved to be
the only codimension one bifurcations.

\begin{figure}[htb]
\centerline{\hbox{
\epsfig{figure=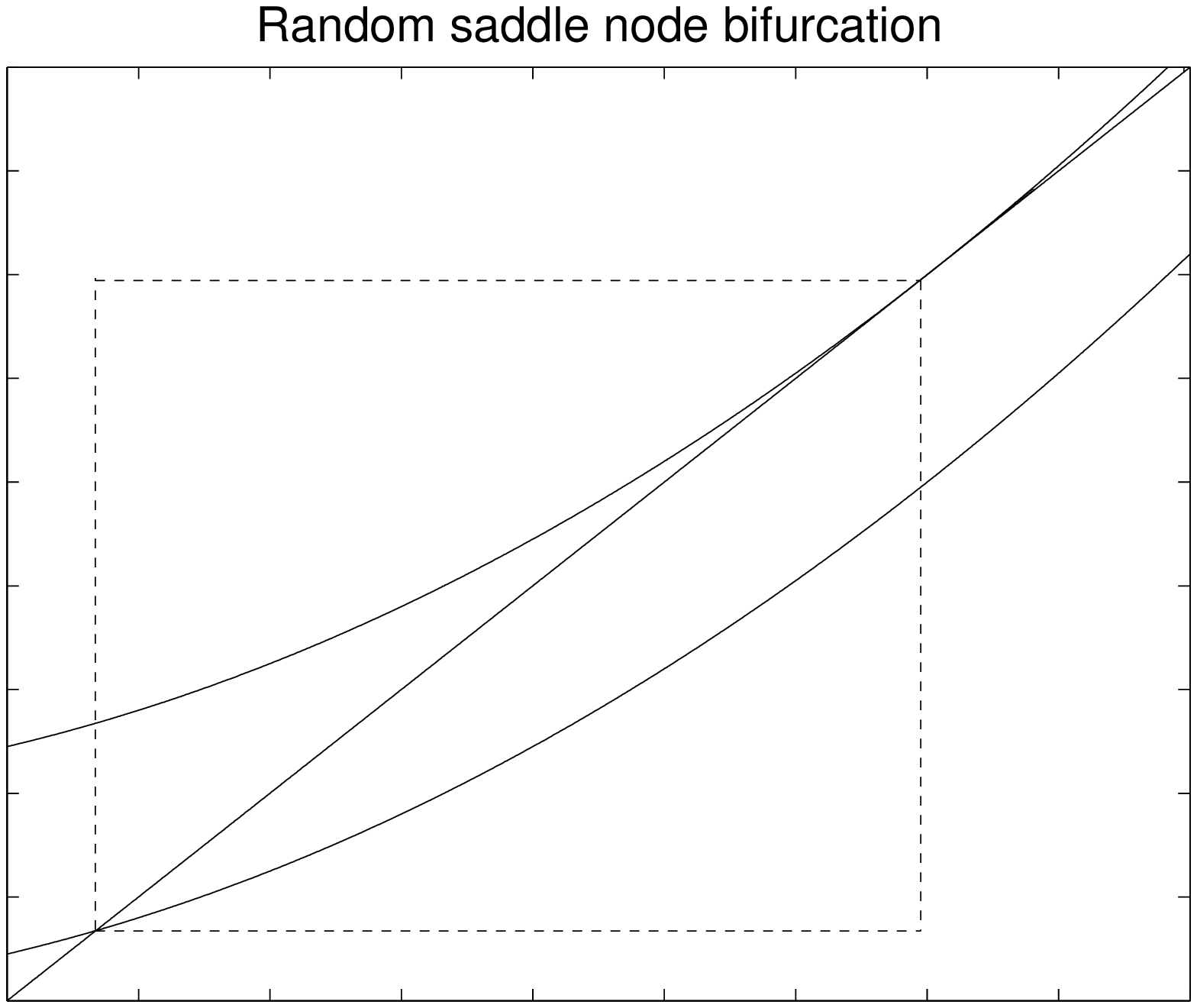,height=8cm,width=8cm}
\epsfig{figure=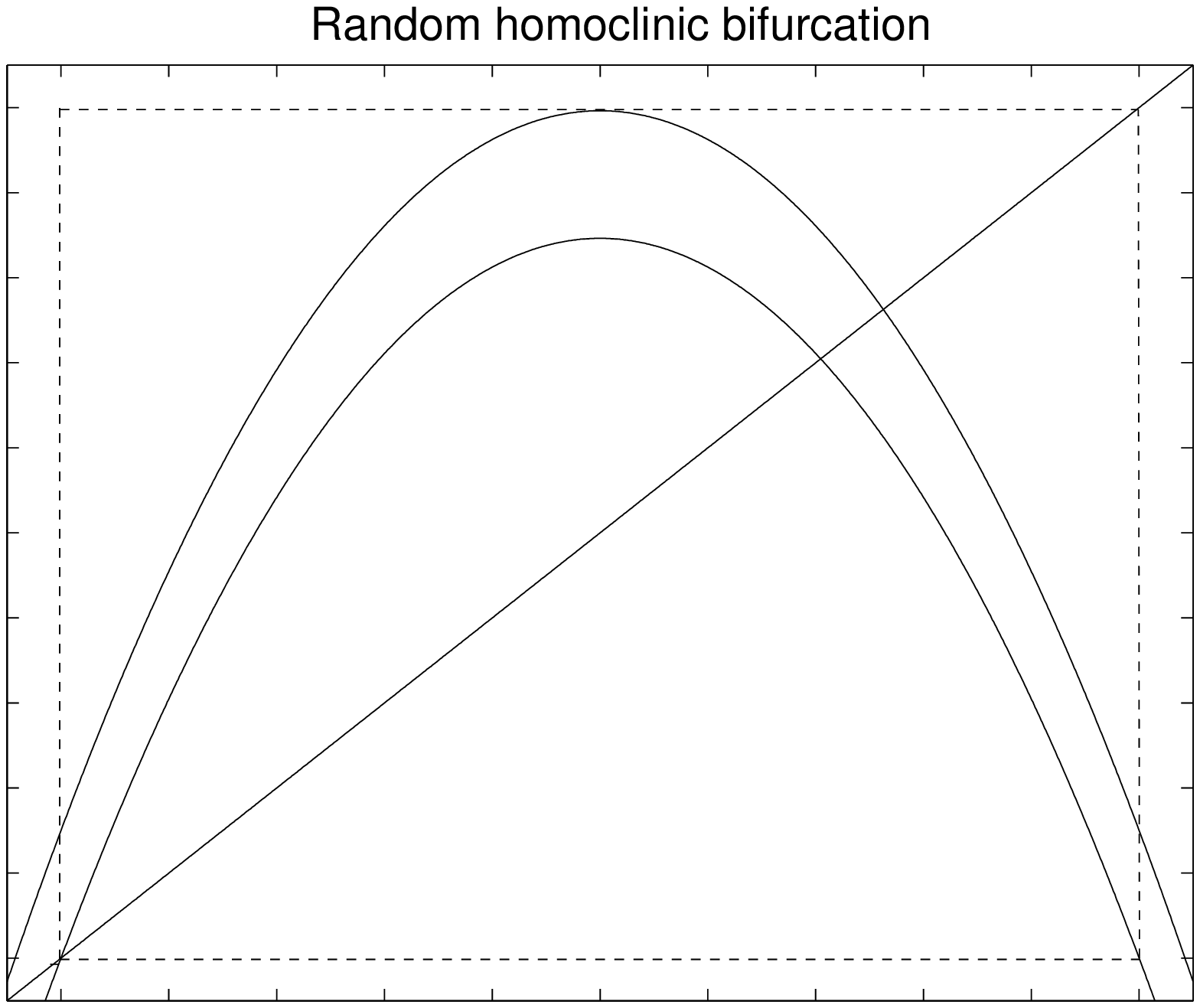,height=8cm,width=8cm}
}}
\caption{\small Consider a random map $f_{a_0}$ for which
points are mapped randomly into the region bounded by the two graphs.
Depicted on the left are the graphs a random map $f(\cdot;\omega)$, 
$\omega \in \partial \Delta$, with a random saddle node bifurcation.
The support of the stationary density
is the interval between the hyperbolic fixed point of the lower map and
the saddle node fixed point of the upper map.
The right picture shows graphs of a random map with a random homoclinic bifurcation.
Here the support of the stationary density stretches from
the left hyperbolic fixed point of the lower map to
the critical value of the upper map.
\label{fig_bifurcations}}
\end{figure}

\begin{definition}\label{def_saddlenode}
The smooth one parameter family of random endomorphisms $f_a$ on the circle
undergoes a {\em random saddle node
bifurcation} at $a=a_0$,
if there exists
$\bar{x}$ in the boundary of the support of
a stationary measure such that
\begin{equation}
f^k_{a_0}(\bar{x};\omega_1,\ldots,\omega_k) = \bar{x},
\qquad
\frac{d}{dx} f^k_{a_0}(\bar{x};\omega_1,\ldots,\omega_k) = 1
\end{equation}
for some $\omega_1,\ldots,\omega_k \in \partial \Delta$.
The random saddle node bifurcation is said to unfold generically, if
\begin{equation}
\left(\frac{d}{dx}\right)^2 f^k_{a_0}(\bar{x};\omega_1,\ldots,\omega_k) \ne 0, \qquad
\frac{\partial}{\partial a} f^k_a (\bar{x};\omega_1,\ldots,\omega_k) \ne 0
\end{equation}
at $a = a_0$.
\end{definition}

\begin{definition}\label{def_homoclinic}
The smooth one parameter family of random endomorphisms $f_a$ on the circle undergoes a
{\em random homoclinic bifurcation} at $a=a_0$,
if there exists
\begin{itemize}
\item a stationary measure $\mu$ with support $E$ with
a hyperbolic periodic point
$\bar{x}_a$ in the boundary of $E$ for all $a$ near $a_0$, and
\item  a critical point
$\bar{y}_a$ for $f_a(\cdot;\omega_1)$, $\omega_1 \in \partial \Delta$,
in the interior of $E$,
\end{itemize}
such that
\begin{equation}
f^l_{a_0} (\bar{y}_{a_0}; \omega_1,\ldots,\omega_l) = \bar{x}_{a_0}
\end{equation}
for some $\omega_2,\ldots,\omega_l \in \partial \Delta$.
The random homoclinic bifurcation unfolds generically if
\begin{equation}
\frac{\partial}{\partial a}
\left(
f^l_{a} (\bar{y}_{a}; \omega_1,\ldots,\omega_l) - \bar{x}_a
\right) \ne 0
\end{equation}
at $a = a_0$.
\end{definition}

\begin{figure}[htb]
\centerline{\hbox{
\epsfig{figure=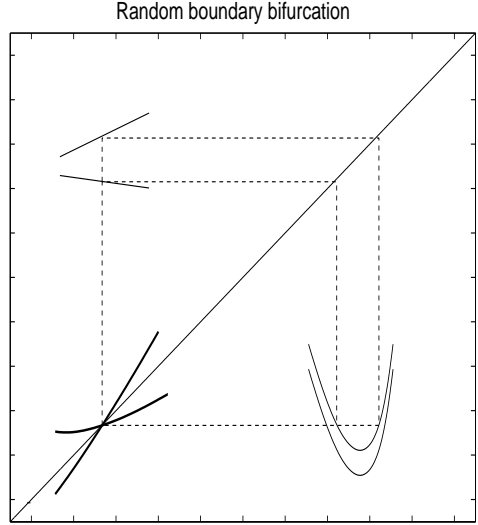,height=8cm,width=8cm}
}}
\caption{\small Depicted are 
parts of the graphs of a random map $f(\cdot;\omega)$, $\omega \in \partial \Delta$. 
The solid curves lie on two of the graphs of 
$f(f(\cdot;\omega_1);\omega_2)$, $\omega_1,\omega_2 \in \partial \Delta$,
intersecting in a point that lies on two hyperbolic
period two orbits (one stable, one unstable) distinguished by different
$\omega$ values. A random boundary bifurcation results 
if this point lies on the boundary of the support of a stationary measure.
\label{fig_boundary}}
\end{figure}

\begin{definition}\label{def_boundary}
The smooth one parameter family of random endomorphisms $f_a$ on the circle
undergoes a {\em random boundary 
bifurcation} at $a=a_0$,
if there exists
$\bar{x}$ in the boundary of the support of
a stationary measure 
and
$(\omega_1,\ldots,\omega_k) \ne (\tilde{\omega}_1,\ldots,\tilde{\omega}_k) 
\in (\partial \Delta)^k$,
such that
\begin{equation}
f^k_{a_0}(\bar{x};\omega_1,\ldots,\omega_k) = \bar{x},
\qquad
\frac{d}{dx} f^k_{a_0}(\bar{x};\omega_1,\ldots,\omega_k) \in (0,1)
\end{equation}
and
\begin{equation}
f^k_{a_0}(\bar{x};\tilde{\omega}_1,\ldots,\tilde{\omega}_k) = \bar{x},
\qquad
\frac{d}{dx} f^k_{a_0}(\bar{x};\tilde{\omega}_1,\ldots,\tilde{\omega}_k) \in (1,\infty)
\end{equation}
Write $\bar{x}_a(\omega_1,\ldots,\omega_k)$
and $\bar{x}_a(\tilde{\omega}_1,\ldots,\tilde{\omega}_k)$
for the continuations of the hyperbolic periodic points.
The random  boundary bifurcation is said to unfold generically, if
\begin{equation}
\frac{\partial}{\partial a} f^k_a (\bar{x}_a(\omega_1,\ldots,\omega_k);\omega_1,\ldots,\omega_k)
\ne \frac{\partial}{\partial a} f^k_a (\bar{x}_a(\tilde{\omega}_1,\ldots,\tilde{\omega}_k;)\tilde{\omega}_1,\ldots,\tilde{\omega}_k)
\end{equation}
at $a = a_0$.
\end{definition}


For an open interval $I$, write $R^k(I,\mathcal{M})$ for the space of $C^k$ families
of random maps in $R^k(\mathcal{M})$ depending on a parameter in $I$.
Equip the space $R^k(I,\mathcal{M})$ with the $C^k$ topology.

\begin{theorem}\label{bifurcation_1D}
For $\{f_a\}$ from an open and dense subset of
$R^\infty(I,{\mathbb S}^1)$,
$f_a$ has only finitely many bifurcations.
A bifurcation point is a
random saddle node bifurcation, a
random homoclinic  bifurcation, or a random boundary bifurcation
and is generically unfolding.
If the number of stationary measures is locally constant at
a bifurcation point, the bifurcation is an intermittency bifurcation.
Otherwise the bifurcation is a transient bifurcation.
\end{theorem}

\noindent {\sc Proof.}
For a bifurcation value for a family in $R^\infty(I,{\mathbb S}^1)$,
either \ref{h_1}, \ref{h_2}, or \ref{h_3} is violated.

Similar transversality arguments as in the proof of Theorem~\ref{stable_1D} show the following. 
For an open and dense subset of
$R^\infty(I,{\mathbb S}^1)$, at most one of these conditions is violated
at a bifurcation value and the resulting bifurcation 
unfolds generically as stated in Definition~\ref{def_saddlenode}, \ref{def_homoclinic} 
or \ref{def_boundary}.
Since the random bifurcations are unfolding generically, they occur isolated.
\qed

\section{Case studies}\label{sec_casestudies}

In this section we illustrate the general theory on two examples;
a randomized version of standard circle diffeomorphisms and
a randomized version of logistic maps on the interval.
We explain how random saddle node bifurcations occur
in both examples and
random homoclinic bifurcations in
random logistic maps.
For the random circle diffeomorphisms we consider rotation numbers
and study their dependence on parameters.
The reader is referred to \cite{melstr93} for the theory of
deterministic circle and interval maps.

\subsection{Random circle diffeomorphisms}\label{sec_circle}

The standard circle map acting on
$x \in {\mathbb R} / {\mathbb Z}$
and depending on parameters $a, \varepsilon$
is given  by
$$
f_a (x) = x + a + \frac{\varepsilon}{2\pi} \sin(2 \pi x) \mod 1.
$$
Consider $f_a$ for a fixed value of $\varepsilon \in (0,1)$ for which
$f_a$ is a diffeomorphism.
%
Introduce the lift $F_a: \mathbb{R}\mapsto\mathbb{R}$,
$$
F_a(x) = x + a + \frac{\varepsilon}{2\pi}\sin(2 \pi x).
$$
It is well known that the rotation number $\rho_a$ of $f_a$,
\begin{equation}\label{rotation1}
\rho_a = \lim_{k \to \infty}\frac{F_a^k(x)- x}{k},
\end{equation}
is well defined and independent of $x$.
The rotation number depends continuously on $a$.
The rotation number is rational precisely if
$f_a$ possesses periodic orbits.
For a fixed rational number $r$, the rotation number of $f_a$
equals $r$ for an interval of $a$ values. In the interior of such an interval,
$f_a$ has exactly one hyperbolic periodic attractor and one hyperbolic
periodic repeller, see \cite{mil00}.

In the following we consider standard circle diffeomorphisms
with a random parameter:
\begin{equation}\label{rarnold}
f_a(x;\omega) = x + \frac{\varepsilon}{2\pi} \sin(2 \pi x) + a + \sigma \omega
\end{equation}
for $x \in {\mathbb R}/{\mathbb Z}$ and
a random parameter $\omega$ chosen from a uniform distribution
on $\Delta = [-1,1]$.
The value of $\sigma$ determines the amplitude of the noise, we
assume it has a fixed value.
We consider fixed $\varepsilon \in (0,1)$ for which $x \mapsto f_a(x;\omega)$ is a diffeomorphism.
Write
\begin{equation}\label{lift1}
F_a(x;\omega)= x + a+\sigma \omega + \frac{\varepsilon}{2\pi}\sin(2\pi x)
\end{equation}
for the lift of $f_a(x;\omega)$.
Note that $F_a (x;\omega) - x$ is periodic in $x$ with period one.

\begin{proposition}
For each parameter value $a$, the random standard circle family
$\{f_a\}$
has a unique stationary measure $\mu_a$.
The density $\phi_a$ of $\mu_a$ is smooth and depends smoothly on $a$.
The support of $\mu_a$ is either the entire circle or
finitely many intervals
strictly contained in the circle.
The latter possibility is only possible
if $\rho_b$ is rational for each $b \in [a-\sigma,a+\sigma]$.
Bifurcations where the support of $\mu_a$ changes discontinuously,
are generic saddle node bifurcations.
There are finitely many such bifurcations.
\end{proposition}

\begin{remark}
Observe that $f_a$ has a hyperbolic fixed point for
$a \in (-\frac{\varepsilon}{2 \pi} ,\frac{\varepsilon}{2 \pi})$.
Hence, $f_a$ has a stationary measure supported on a single interval
precisely if both $a - \sigma > -\frac{\varepsilon}{2 \pi}$
and $a + \sigma < \frac{\varepsilon}{2 \pi}$. This occurs for a nonempty
interval of $a$ values if $\sigma < \frac{\varepsilon}{2 \pi}$.
\end{remark}

\noindent {\sc Proof.}
It is well known that a circle diffeomorphism with irrational rotation number has its orbits
lying dense in ${\mathbb R}/{\mathbb Z}$.
It follows that
if the family of circle maps $f_a (x ;\omega)$ for varying $\omega \in \Delta$
contains a member with irrational rotation number,
there is a (necessarily unique) stationary measure supported on all of
${\mathbb R}/{\mathbb Z}$.

Suppose now that $f_a (x ;\omega)$ for  each $\omega \in \Delta$
has rational rotation number $\rho_{a + \sigma \omega} = p/q$.
Write $x_{\omega}$ for a periodic point from a periodic attractor
of  $x\mapsto f_a (x ;\omega)$ depending continuously on $\omega$.
Recall that $x\mapsto f_a (x ;\omega)$ has a unique periodic
attractor. Let $V_a=\cup_{\omega \in \Delta} x_{\omega}=[x_{-1},
x_{1}]$. The random standard family is increasing in $x$ and in
$\omega$, so that for all $x \in V_a$, and all ${\boldsymbol
\omega} \in \Delta^{\mathbb N}$ we have
$$
x_{-1}= F_a^q(x_{-1}; -1) \leq F_a^q(x_{-1};{\boldsymbol
\omega})\leq F_a^q(x;{\boldsymbol \omega}) \leq
F_a^q(x_{1};{\boldsymbol \omega}) \leq F_a^q(x_{1};1)= x_{1}.
$$
It follow that the orbit of $V_a$ is invariant.
For a fixed $\omega \in \Delta$,
all points outside the unique periodic repeller of $f_a(\cdot; \omega)$
are attracted to its periodic attractor.
This implies that there is a unique
stationary measure supported on the orbit of $V_a$.

Compute
$$
\frac{\partial}{\partial a} f_a^k (x ; {\boldsymbol \omega}) =
\sum_{i=0}^k \frac{\partial}{\partial a}f (f^i(x;{\boldsymbol
\omega};{\boldsymbol \omega}))
              \frac{d}{dx}f^i_a (f^{k-i}_a (x;{bf \omega};{\boldsymbol \omega})).
$$
As all terms in the sum are positive, a random saddle node bifurcation
occurs isolated. The random family $\{ f_a\}$ therefore has only
a finite number of random saddle node bifurcations.
\qed\\

\begin{figure}[htb]
\centerline{\hbox{ \epsfig{figure=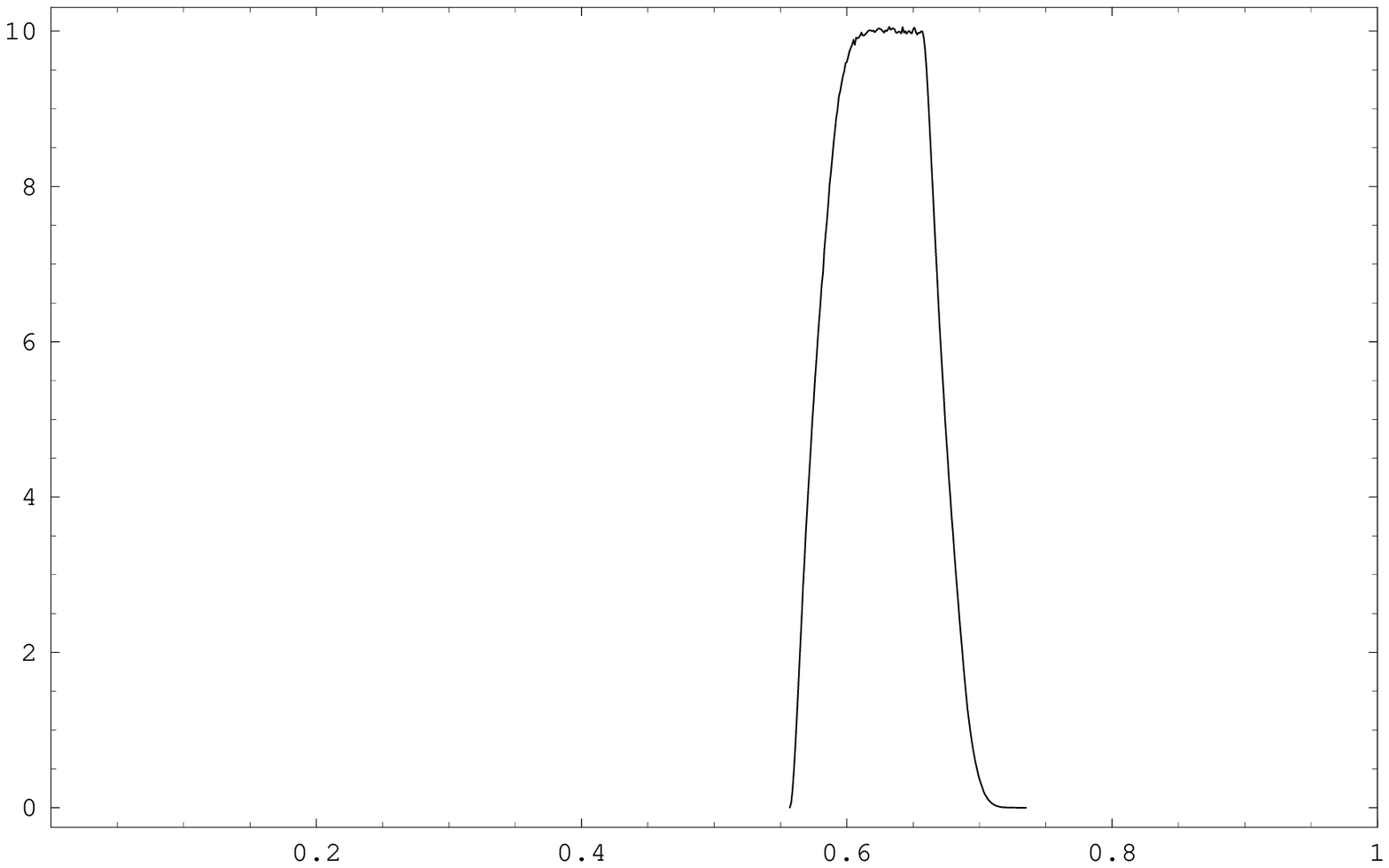,height=6cm,width=8cm}
\epsfig{figure=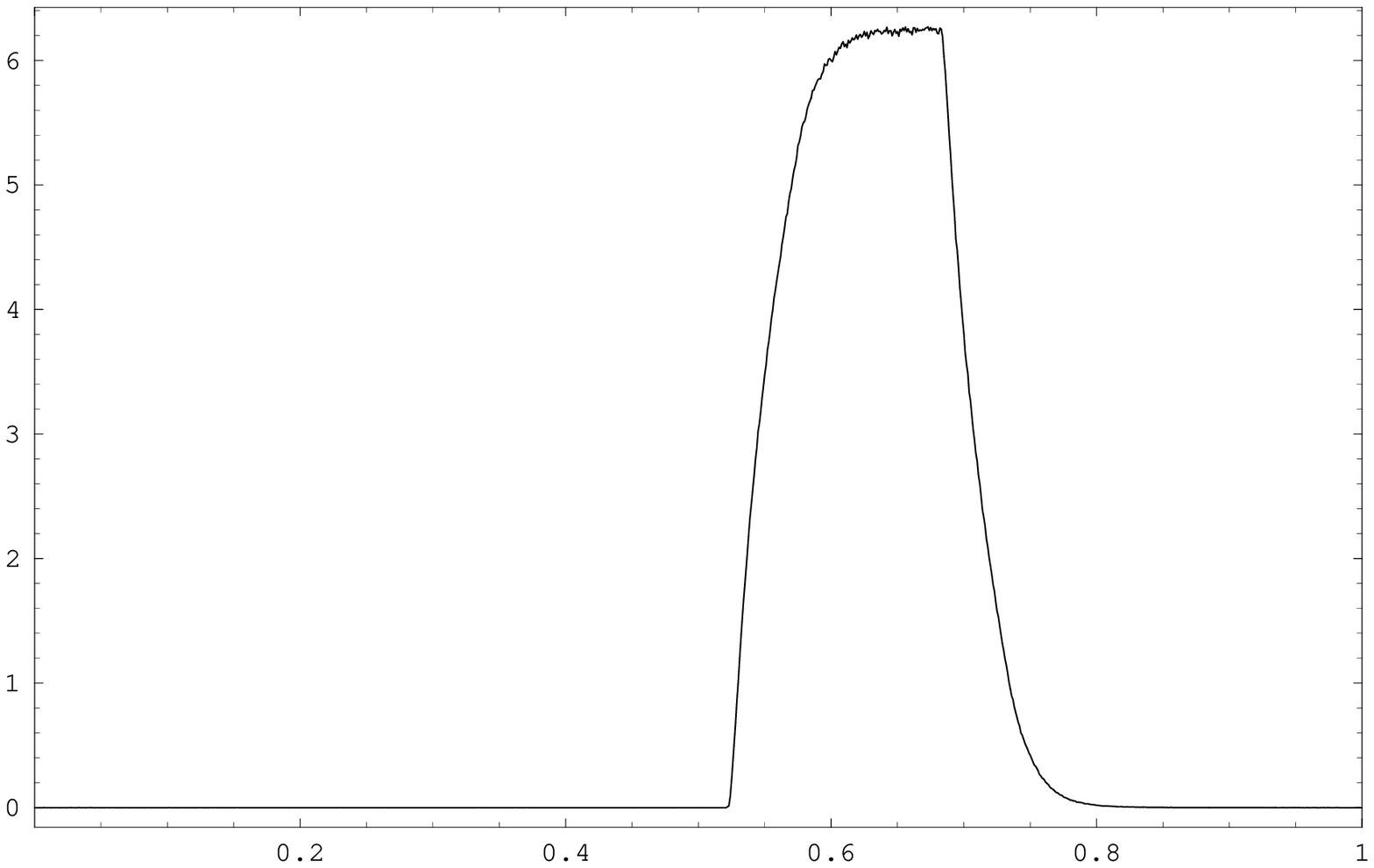,height=6cm,width=8cm}}} \caption{\small
Numerically computed stationary densities of the 
random standard circle map.
On the left for $|a|+\sigma <\varepsilon/2\pi$, on
the right for $|a|+\sigma > \varepsilon/2\pi$.  
The explosion of the
support of the stationary density follows a random saddle node bifurcation. \label{circ_bif}}
\end{figure}

We define the rotation number for the random standard circle map,
when its exist, by
\begin{equation}\label{rrotation}
 \rho_a(x;{\boldsymbol \omega})  =  \lim_{k \to \infty}\frac{F_a^k(x ;{\boldsymbol \omega})-x}{k}.
\end{equation}
The rotation number measures the average rotation per iterate of $f_a$.
Note that $\rho_a$ is a random variable,
depending also on the starting point $x$.

A simple but usefull lemma shows that $\rho_a$ is independent
of the initial condition $x$.

\begin{lemma}\label{Initcond}
If $\rho_a (x; {\boldsymbol \omega})$ exists for some $x \in
{\mathbb R}/{\mathbb Z}, {\boldsymbol \omega} \in \Delta^{\mathbb
N}$, then $x \mapsto \rho_a(x;{\boldsymbol \omega}$ exists for all
$x \in {\mathbb R}/{\mathbb Z}$ and is constant in $x$.
\end{lemma}

\noindent {\sc Proof.} Observe that $F_a^k( \cdot;\omega)$ is a
lift of $f_a^k(\cdot;{\boldsymbol \omega})$, so that $F_a^k(
x;\omega) - x$ is periodic in $x$ with period 1. Thus
$$
\max_{x \in \mathbb{R}}\{F_a^k( x;\omega)- x\}-\min_{x \in \mathbb{R}}\{F_a^k(x ;\omega) - x\}<1.
$$
Compute
\begin{align*}
|F_a^k(x ;{\boldsymbol \omega})-F_a^k( y;{\boldsymbol \omega})|& \leq |(F_a^k(x ;{\boldsymbol \omega})-x)-(F_a^k( y;{\boldsymbol \omega})-y)|+|x-y|\\
 &\leq 1 + |x-y|,
\end{align*}
so that
$$
\lim_{k \to \infty}\Bigl(\frac{F_a^k(x ;{\boldsymbol \omega}) -
x}{k}
                         -\frac{F_a^k(y ;{\boldsymbol \omega}) - y}{k}\Bigr)=0.
$$
It follows that
the limit $\lim_{k \to \infty}\frac{F_a^k(x ;\omega) - x}{k}$, if it exists,  is independent of $x$.
\qed\\

Write
$$
F_a(x;{\boldsymbol \omega})=x+\delta_a(x;{\boldsymbol \omega}),
$$
where the function $\delta_a(x;{\boldsymbol \omega})$ is periodic
with period one in the variable $x$. We can consider $\delta$ as a
function  defined on $\mathbb{R}/\mathbb{Z}$. A simple induction
argument gives for each $k \in \mathbb{N}$,
\begin{equation}\label{birkhoff}
f_a^k(x;{\boldsymbol \omega})= x + \sum_{i=0}^{k-1}\delta \circ
S^i(x;{\boldsymbol \omega})
\end{equation}
where $S$ is the skew product system (see equation (\ref{skew})) on
$\mathbb{R}/\mathbb{Z} \times \Delta^{\mathbb N}$ .
Recall that $\mu_a \times \nu^{\infty}$ is an $S$-invariant measure.
\begin{figure}[htb]
\centerline{\hbox{ \epsfig{figure=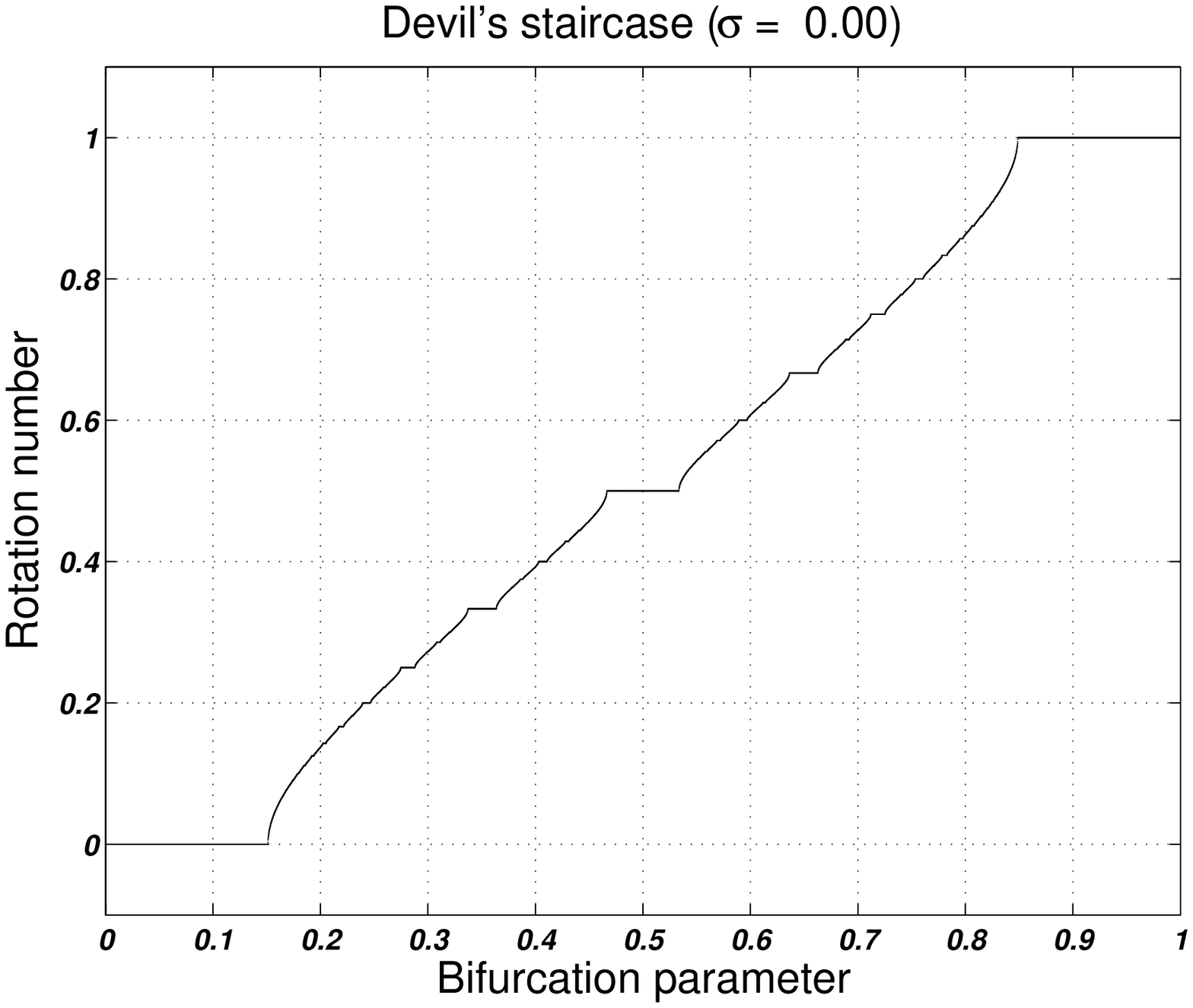,height=6cm,width=8cm}
\epsfig{figure=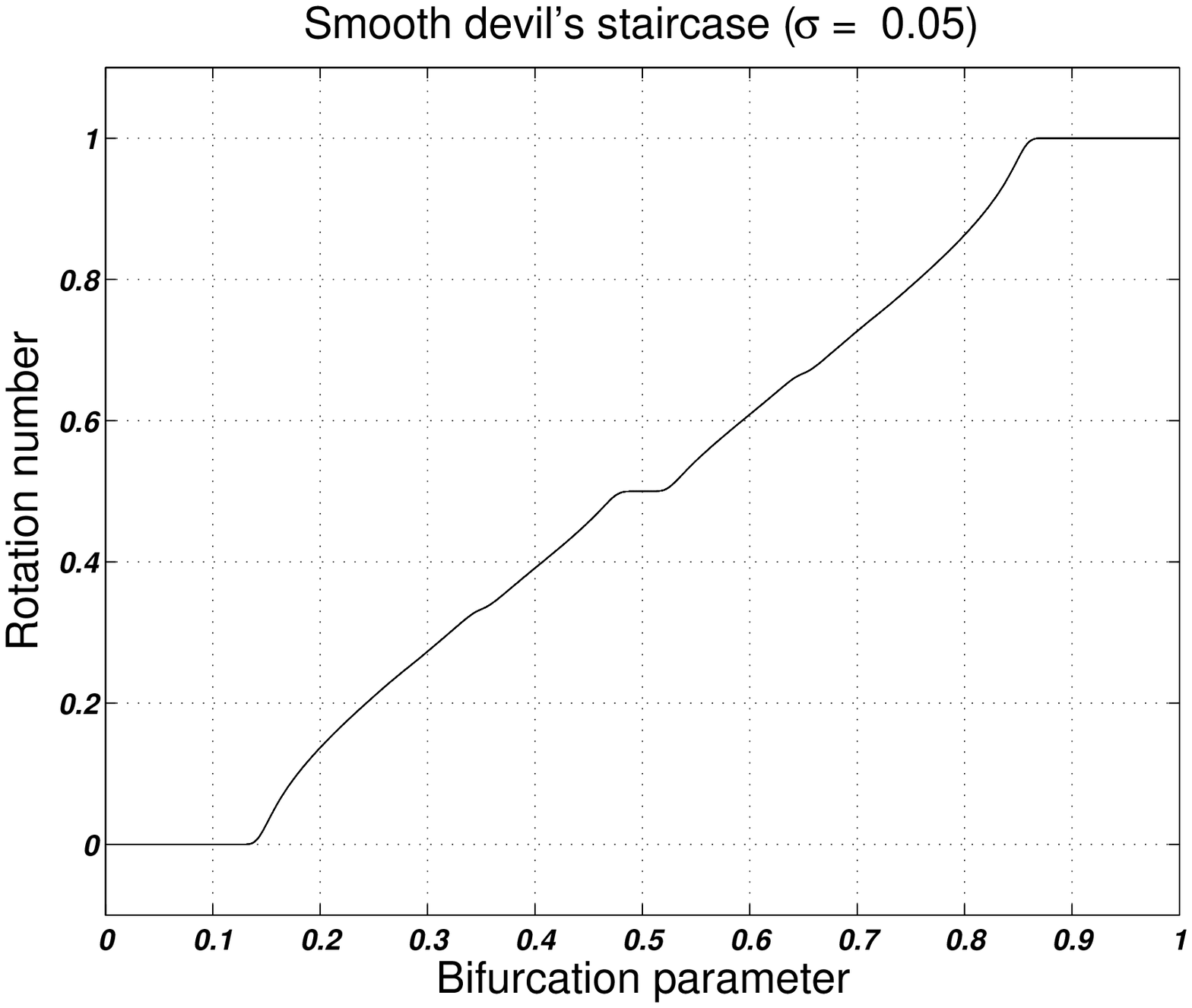,height=6cm,width=8cm} }}
\caption{\small The function $a\mapsto\rho_a$. On the
left the devil's staircase; the rotation number
of the deterministic standard family.
On the right
the rotation number of the random standard family. \label{fig_Devil}}
\end{figure}

\begin{proposition}
\begin{equation}\label{rotation2}
\rho_a = \int_{\mathbb{R}/\mathbb{Z}}\mathbb{E}(\delta(x;\omega))d\mu_{a}(s) \qquad \nu^{\infty}-a.s.
\end{equation}
where $\mathbb{E}$ is the expectation operator.
The right hand side of (\ref{rotation2}) is independent of $x$ and is
a smooth and nondecreasing function of $a$.
\end{proposition}

\noindent {\sc Proof.} Consider the Birkhoff sum in equation (\ref{birkhoff})
$$
\frac{f_a^k(x;\omega)- x}{k}=\frac{1}{k}\sum_{i=0}^{k-1}\delta_a \circ S^i(s;\omega).
$$
By Birkhoff's ergodic theorem,
\begin{align}
\nonumber \lim_{k \to \infty}\frac{f_a^k(x;\omega)- x}{k}& =\int_{ \mathbb{R}/\mathbb{Z}}\int_{\Omega}
\beta_a(s;\omega)\mu_a(ds)\nu^{\infty}(d\omega) \qquad \mu_a \times \nu^{\infty}-a.s.\\
\label{rotation3} &=\int_{ \mathbb{R}/\mathbb{Z}}\mathbb{E}(\beta_a(s;\omega))\mu_{a}(ds)
                              \qquad \mu_a \times \nu^{\infty}-a.s.
\end{align}
The fact that the rotation number when it exist is independent of the initial
point, see (\ref{Initcond}), implies that this equality holds for all $x$ and $\nu^{\infty}-a.s.$.

Write
$a \mapsto h(a)$ for the right hand side of (\ref{rotation2}).
Smoothness of $h$ follows from smoothness of the stationary density $\phi_a$
and (\ref{rotation3}).
Write the standard family as $R_a\circ f(.;\omega)$ where $f(.;\omega)$ is the random map
$f(x;\omega)=x+ \frac{\varepsilon}{2\pi}\sin(2\pi x)+\sigma\xi(\omega)$ and $R_a$ is the translation with
coefficient $a$.
Then for $a_1<a_2$ and $k \geq 1$,
$(R_{a_1}\circ f(.; \omega))^k< (R_{a_2}\circ f(.;\omega))^k$ and thus,
$$
\rho_{a_1} =\lim_{k \to \infty}\frac{(R_{a_1}\circ f(.;\omega))^k-Id}{k}\leq\lim_{k \to
\infty}\frac{(R_{a_2}\circ f(.; \omega))^k-Id}{k}=\rho_{a_2}.
$$
\qed

\subsection{Random unimodal maps}\label{sec_logistic}

This section is devoted to the investigation of the randomized version of the logistic family
\begin{equation}
f_a (x;\omega) = (a + \sigma \omega) x (1-x).
\end{equation}
on $[0,1]$.
The  random parameter $\omega$ be chosen from a uniform
distribution on $\Delta = [-1,1]$.
Throughout this section we will assume that
\begin{equation}\label{standingassumption}
a + \sigma \omega \in (1,4),
\end{equation}
for all $\omega \in \Delta$.
As a consequence, the interval $[0,1]$ is mapped
into itself by each map $x\mapsto f_a(x;\omega)$
and the fixed point at the origin is repelling.
Any stationary measure will therefore have support contained in $(0,1)$.
We will demonstrate that there is only one stationary measure.

\begin{proposition}
The random logistic map $f_a$ has a unique stationary measure.
\end{proposition}

\noindent {\sc Proof.}
We collect some facts from unimodal dynamics needed in the sequel of the  proof.
The following facts hold for unimodal maps with negative Schwarzian
derivative such as the logistic map $x \mapsto a x (1-x)$.
By Guckenheimer's theorem, see \cite[Theorem~III.4.1]{melstr93},
$x \mapsto a x (1-x)$ possesses  a unique
attractor $\Lambda_a$.
The attractor $\Lambda_a$ is either a periodic  attractor, a solenoidal attractor,
or a finite union of intervals on which the map acts transitively.
In all cases, the omega-limit set of the critical point $c$  (with $c= \frac{1}{2}$
for the logistic map)
is contained in $\Lambda_a$. In fact, if $\Lambda_a$ is not a periodic attractor,
then $c$ is contained in $\Lambda_a$.
It follows from a result of Misiurewicz, see \cite[Theorem~III.3.2]{melstr93},
that the basin of attraction of $\Lambda_a$ is an open and dense subset of $(0,1)$.

We will distinguish the following two cases.
\begin{description}
\item[Case (i):] There exists $\omega \in \Delta$, so that $c \in \Lambda_{a + \sigma \omega}$,
\item[Case (ii):] otherwise.
\end{description}
The two cases are treated separately.

\noindent {\em Case (i)}:
Write
$$
W =  \bigcap_{n \ge 0} \overline{\bigcup_{i \ge n} f^i_a(c;\Delta^{\mathbb N})}
$$
for the omega-limit set of $c$ under all possible random
iterations. Observe that $W$ is an invariant set. From the
properties of the noise, $W$ consists of a finite union of
intervals. Note that $c \in W$, so that $W$ equals the closure of
the positive orbit $\overline{\bigcup_{i \ge 0}
f^i_a(c;\Delta^{\mathbb N})}$ of $c$ under all possible random
iterations. We will prove that for each $x \in (0,1)$, $y \in W$
and $\varepsilon >0$, there exist $n >0$ and $\bar{{\boldsymbol
\omega}} \in \Delta^{\mathbb N}$ with the property that
\begin{equation}\label{orbits}
|f^n_a ( x ; \bar{{\boldsymbol \omega}}) - y | < \varepsilon.
\end{equation}
This implies the $W$ is the unique minimal invariant set, which in turn implies
the theorem in the first case.

Fix $x \in (0,1)$, $y\in W$, $\varepsilon>0$. From the
construction of $W$, there exist ${\boldsymbol \omega}_1 \in
\Delta^{\mathbb N}$, $i>0$, so that $| f^i_a( c ; {\boldsymbol
\omega}_1) - y | < \varepsilon$. By continuity of $x \mapsto
f^i_a(x;{\boldsymbol \omega}_1)$, the same holds with $c$ replaced
by a point from a $\delta$ neighborhood of $c$ for some $\delta
>0$. We need to establish the existence of $\hat{{\boldsymbol \omega}} \in
\Delta^{\mathbb N}$ and $j>0$ so that $| f^j_a(x;\hat{{\boldsymbol
\omega}}) - c| < \delta$. Let $\omega_2 \in \Delta$ be such that
$c \in \Lambda_{a + \sigma \omega_2}$. Since the basin of
attraction of $ \Lambda_{a + \sigma \omega_2}$ is open and dense,
there exists $\omega_3 \in \Delta$ with $x_1 = f_a(x ; \omega_3)$
contained in the basin of attraction of $ \Lambda_{a + \sigma
\omega_2}$. For $i$ large, $x_i = f^{i-1}_a(x_1;\omega_2,
\omega_2,\cdots)$ is as close as desired to $\Lambda_{a + \sigma
\omega_2}$. If $\Lambda_{a + \sigma \omega_2}$ is a finite union
of intervals, we get that $x_i$ is contained in $\Lambda_{a +
\sigma \omega_2}$ for large enough $i$.
As inverse images of $c$ for $x \mapsto f_a(x;\omega_2)$ are dense
in $\Lambda_{a + \sigma \omega_2}$, one deduces that there exist
${\boldsymbol \omega}_4 \in \Delta^{\mathbb N}$ and $k >0$ so that
$f^k_a(x_{i+1}; {\boldsymbol \omega}_4)$ lies in a $\delta$
neighborhood of $c$. Indeed, if $\Lambda_{a + \sigma \omega_2}$ is
a finite union of intervals, then we find $\omega_5$ with $x_{i+1}
= f(x_i;\omega_5)$ equal to an inverse image of $c$. If
$\Lambda_{a + \sigma \omega_2}$ is a solenoidal attractor, then $x
\mapsto f_a(x;\omega_2)$ is infinitely renormalizable. In this
case one can use ${\boldsymbol \omega}_4 = (\omega_2,
\omega_2,\cdots)$. Also for a periodic attractor $\Lambda_{a +
\sigma \omega_2}$ containing $c$ one uses ${\boldsymbol \omega}_4
= (\omega_2, \omega_2,\cdots)$.

\noindent {\em Case (ii)}:
By Guckenheimer's theorem,
$x \mapsto f_a(x;\omega)$ possesses a unique periodic attractor
for each $\omega \in \Delta$.
Write $V$ for the
the union of $\Lambda_{a + \sigma \omega}$
over $\omega \in \Delta$.
Define
$$
W = \overline{\bigcup_{i\ge 0} f^i_a(V;\Delta^{\mathbb N})}.
$$
This is clearly an invariant set.
Note that we do not claim that $c$ is outside of $W$.
Arguments as before prove (\ref{orbits}) with this definition of $W$: for suitable noise one
finds an orbit starting at $x$ that approaches a point in $V$ and then with further
iterates approaches $y \in W$.
\qed

\begin{remark}
The above proof applies to show that
a random unimodal map $g(x;\omega)$ with negative Schwarzian derivative
for each $\omega$
(the invoked theorem by Guckenheimer is true for these maps) has a unique
stationary measure.
\end{remark}


\begin{figure}[htb]
\centerline{\hbox{
\epsfig{figure=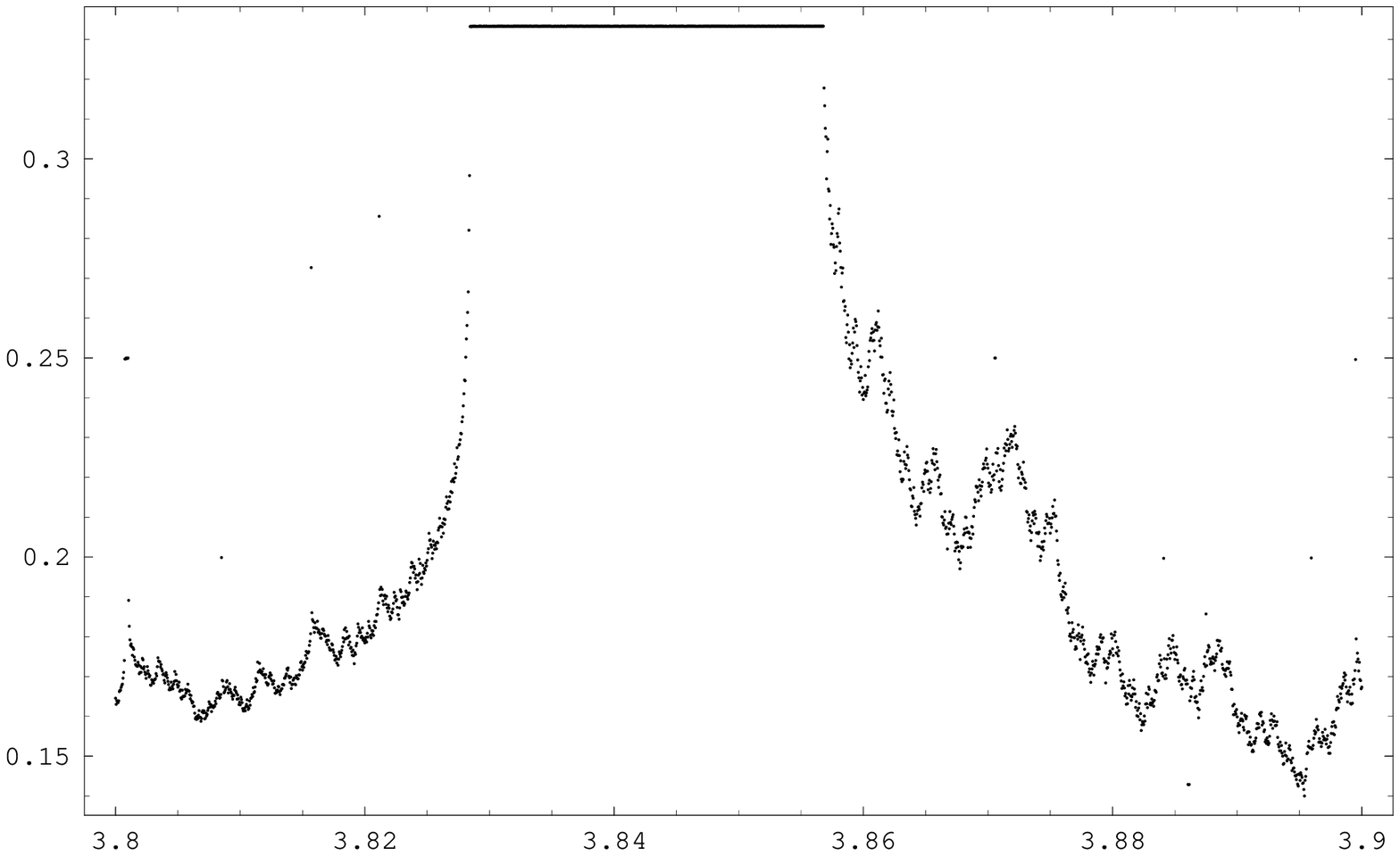,height=6cm,width=7.5cm}
\epsfig{figure=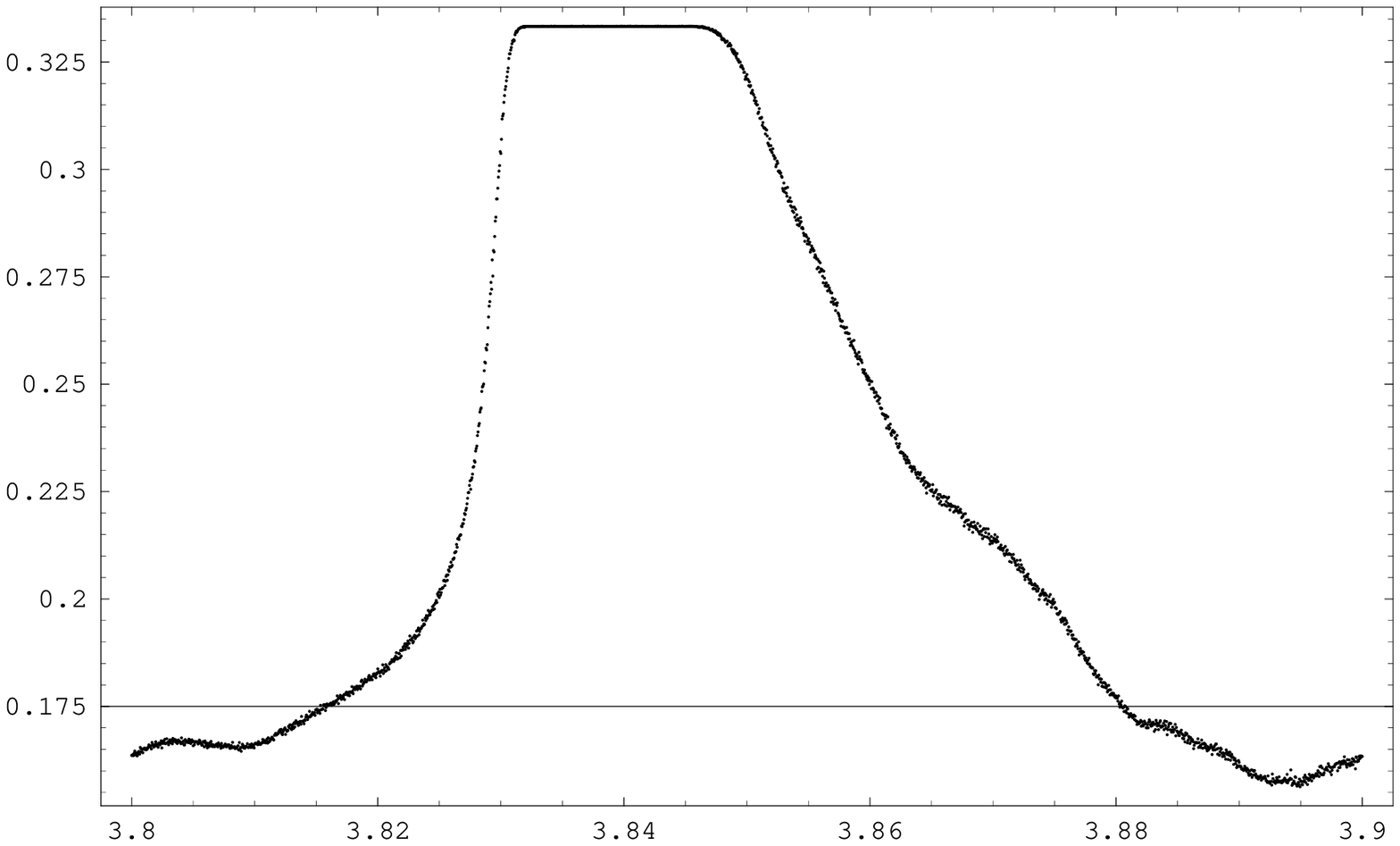,height=6cm,width=7.5cm}
}}
\caption{\small 
Numerically computed Birkhoff averages 
$\lim_{n\to \infty} \frac{1}{n} 
    \sum_{i=0}^{n-1} \phi (f_a^i(x;{\boldsymbol \omega}))$
of $\phi(x) = 1_{[0.4,0.6]}$ for the
logistic family (left picture) and the 
random logistic family with $\sigma=0.005$ (right picture), 
for parameters $a$ ranging from 3.8 to 3.9. The flat part in the left picture, 
where the average equals 1/3,
runs from a saddle node bifurcation to a homoclinic bifurcation.
The numerical computations show that 
for the random logistic family these are replaced by their random versions;
for parameters from the flat part the
stationary measure is supported on three disjoint intervals cycled by the
random map.  
\label{fig_ba}}
\end{figure}
Perturbing away from the deterministic logistic family one sees that 
both random saddle node bifurcations and random homoclinic bifurcations occur in the random logistic
family $\{f_a\}$ for small noise levels. 
Typically one can expect the following scenario.
We start by recalling 
some facts concerning the dynamics of the deterministic map $f_a(\cdot;0)$.
The map $f_a (\cdot;0)$ is called renormalizable if there exists an interval $I$ and a positive integer $q$, so that
$f_a^q (I;0) \subset I$.
Let $[a_-.a_+]$ be a maximal interval 
so that $f_a(\cdot;0)$ is renormalizable for $a \in [a_-,a_+]$ 
with $q$ constant. 
Then $f_{a}$ undergoes a saddle node bifurcation at $a=a_-$ involving a periodic orbit of period $q$.
At $a=a_+$, $f_a$ undergoes a homoclinic bifurcation, where an iterate of $f_a$ maps the critical point
onto a periodic orbit of period $q$.
For small noise levels (i.e. $\sigma$ small) one expects a random saddle node bifurcation near $a = a_-$
and a random homoclinic bifurcation near $a=a_+$.
Figure~\ref{fig_ba} illustrates this by computing Birkhoff averages for the logistic family and
the random logistic family. See \cite{homyou02} for explanations of the computations
for the logistic family.

\appendix

\section{Representations of discrete Markov processes}\label{sec_representation}

In this appendix we explore the relation between random maps and
discrete Markov processes given by stochastic transition functions.
The random maps considered in this paper depend
on random parameters, where the number of random parameters
equals the dimension $n$ of the state space $\mathcal{M}$.
Proposition~\ref{representation}
gives a wide class of Markov processes that can be represented by random maps
by $n$ random parameters. The Markov process given by random maps
depending on a larger number of random parameters
(or even given by some measure on the space of maps)
can be represented by
random maps with $n$ random parameters.

Iterating a random map involves more random parameters
obtained by independent draws at each iterate.
By means of an example we explain how random maps with a smaller
number of random parameters may be brought into
the context of this paper.
Consider the delayed
logistic map $x_{n+1} = \mu x_n (1 - x_{n-1})$.
Let $y_{n+1} = x_{n}$. This defines a dynamical system
$(x_{n+1},y_{n+1}) = (\mu x_n (1-y_n) , x_n)$.
Assume now that $\mu$ is a random parameter varying in some interval
with some distribution.
This yields a random map
$$
f(x,y;\mu) = (\mu x (1-y) , x).
$$
The derivative $Df$ is singular along $x=0$.
As $\mu$ is a single random parameter, this
random diffeomorphism does not fit into the context considered in this paper.
Considering two iterates gives two independent draws $(\mu,\nu)$
of the random parameter (that is, random parameters taken from a square)
and yields the random map
$$
f^2(x,y;\mu,\nu) = ( \mu x (1-y) , \nu \mu x (1-x) (1-y) ).
$$
If $x$ and $y$ stay away from $0$ and $1$, the map and the dependence of $(\mu,\nu)$
are injective. The second iterate of the delayed logistic map with
bounded parametric noise fulfills the assumptions in this paper.

\begin{figure}[htb]
\centerline{\hbox{
\epsfig{figure=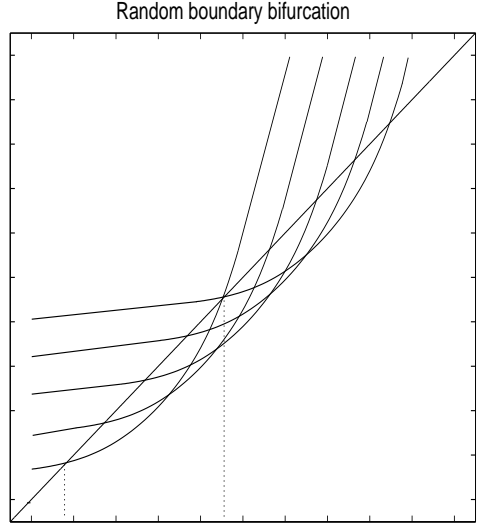,height=8cm,width=8cm}
}}
\caption{\small A random map $f(x;\omega) = f(x-\omega;0) + \omega$
with a random boundary bifurcation of the stationary measure with support between
the ordinates indicated by dotted lines.
If $\omega$ is chosen from a uniform distribution, then the density of 
the transition function will not be bounded.
\label{fig_boundary_2}}
\end{figure}

There are other examples of maps with
parametric noise that cannot be made to fulfill the assumptions
used in this paper. For instance, random maps $f(x;\omega) =  x + (x-\omega)^2$
with random $\omega$ from an interval, fail to satisfy the injectivity assumption
of $\omega \mapsto f(x;\omega)$. If $\omega$ is chosen from a uniform distribution,
then the density of the transition function will not be bounded.
Figure~\ref{fig_boundary_2} indicates a random boundary bifurcation 
for a similar random map.

Consider discrete Markov processes given by transition functions $P(x,\cdot)$.
The following properties hold.
\begin{description}
\item
For fixed $A \in {\cal B}$, $x\mapsto P(x,A)$ is measurable.
\item
For fixed $x \in {\mathbb R}^n$, $P(x,\cdot)$ is a probability measure.
\end{description}
Denote by $y \mapsto k(x,y)$ the density of $P(x,\cdot)$.
Write $U_x$ for the support of $k(x,\cdot)$ and
let $U = \displaystyle{\cup_x} \left(  \{x\}\times U_x \right)$.
Assume
that $U_x$
is diffeomorphic to the closed unit ball $\Delta$ in ${\mathbb R}^n$
and varies smoothly with $x$.
We will assume that $y\mapsto k(x,y)$ depends smoothly on $(x,y) \in U$,
meaning that $k$ can be extended to a smooth
function defined on an open neighborhood of $U$.
Under these conditions we will construct
a representation by a finitely parameterized family of endomorphisms.
That is, we will construct a family of endomorphisms
$\{f_\mu\}$ on  ${\mathbb R}^n$, with parameters
$\mu$ from an $n$ dimensional ball,
and a measure $\nu$ on the parameter space
so that
$P(x,A)$ equals $\nu \{ \mu \in \Delta  \;|\; f_\mu(x) \in A \}$.
A corresponding result holds for discrete Markov processes
with noise from an $n$-dimensional box,
see \cite[Appendix~D]{bondiavia05}.
See \cite{kif86} for a discussion of the existence of
representations by sets of measurable or continuous maps.
The paper \cite{qua91} contains a result on representations by
endomorphisms, under the assumption of unbounded noise.

\begin{proposition}\label{representation}
There is a family of endomorphisms $f_\mu$, $\mu \in \Delta$,
and a measure $\nu$ on $ \Delta$ with smooth strictly positive density,
so that
\begin{enumerate}
\item
$(x,\mu) \mapsto f_\mu(x)$ is smooth,
\item
for each $x \in \mathcal{M}$, $\mu \mapsto f_\mu (x)$ is injective,
\item
$P(x,A) = \nu (\mu \in \Delta \;|\; f_\mu(x) \in A)$.
\end{enumerate}
\end{proposition}

\noindent {\sc Proof.}
We follow the arguments in \cite{bondiavia05}, combined with the use of
polar coordinates to map the unit ball $\Delta$  to $[0,1]^n$.
Let $\psi_x :  V_x \to \Delta$ be a diffeomorphism, depending smoothly on $x$,
from the support $V_x$ of $y \mapsto k(x,y)$ to the unit ball.
Consider polar coordinates $\chi: [0,1]^n \to \Delta$ on the unit ball,
$$
\chi( \xi_1,\ldots,\xi_n) = \xi_1 \left( \begin{array}{c}
 \cos (\pi \xi_2)
\\
 \sin (\pi \xi_2) \cos(\pi \xi_3)
\\
 \vdots
\\
\sin(\pi \xi_2) \cdots \sin(\pi \xi_{n-1} ) \cos(2 \pi \xi_n)
\\
 \sin(\pi \xi_2) \cdots \sin(\pi \xi_{n-1} ) \sin(2 \pi \xi_n) \end{array} \right).
$$

For $\xi \in [0,1]^n$, define sets $B_i (\xi)$, $0\le i \le n$, by
$$
B_i (\xi) = \prod_{j=1}^i [0,\xi_j] \times \prod_{j=i+1}^n [0,1].
$$
Write $C_i(\xi) = \psi_x^{-1} \chi^{-1} (B_i(\xi))$ and let
$\omega = (\omega_1,\ldots,\omega_n)$ be given by
$$
\omega_i = \int_{C_i(\xi)} k(x,y) dm(y) \left/
             \int_{C_{i-1}(\xi)} k(x,y) dm(y)\right..
$$
Since $k > 0$, $\omega = \Theta (\xi)$ gives a 1-1 correspondence.
Let $\eta_i =  \int_{C_i(\xi)} dm(y) \left/
             \int_{C_{i-1}(\xi)}  dm(y)\right.$.
Here $\eta = \Psi(\xi)$ is  a 1-1 correspondence.
The correspondence $\omega  \to \eta$ is a smooth diffeomorphism
as $k$ is smooth and strictly positive.
Then
$$
f_\mu (x) = \psi_x \chi \Psi^{-1} \Theta \chi^{-1} (\mu).
$$
gives the required smooth random maps.
\qed\\

For discrete Markov processes on a circle there is an easy necessary and sufficient
condition on the transition maps for a representation by random
diffeomorphisms.

\begin{proposition}
Let $\mathcal{M}$ be the circle endowed with Lebesgue measure. 
Write $V_x = [l_{-}(x),l_{+}(x)]$.
There is a representation by random smooth diffeomorphisms
if and only if
$$
- k(x,l_{-}(x)) l_{-}'(x)
 +
\int_{l_{-}(x)}^{z} \frac{\partial}{\partial x}  k(x,y) dy \ne 0
$$
for $z \in V_x$.
\end{proposition}

\noindent {\sc Proof.}
The construction of the representation by random smooth
maps proceeds as follows.
For $\xi \in [0,1]$, write
$C  (\xi) = [l_{-}(x) , l_{-}(x) + \xi (l_{+}(x) - l_{-}(x))]$
and let
$$
\omega = \int_{C (\xi)} k(x,y) dy.
$$
Since $k>0$, the map $\Theta$, $\Theta(\xi) = \omega$,  is a diffeomorphism.
The representation by random diffeomorphisms is given through
$$
f_\omega (x) = l_{-}(x) + \Theta^{-1} (\omega) (l_{+}(x) - l_{-}(x)).
$$
Note that for fixed $\omega$,
\begin{equation}\label{propertyomega}
\frac{d}{dx} \int_{[l_{-}(x),f_\omega (x)]}   k(x,y) dy = 0.
\end{equation}
With, say, $l_{-} < f_\omega$, (\ref{propertyomega}) yields
$$
- k(x,l_{-}(x)) l_{-}'(x) + k(x,f_\omega(x)) f_\omega' (x) +
\int_{l_{-}(x)}^{f_\omega (x)} \frac{\partial}{\partial x}  k(x,y) dy = 0.
$$
Hence $f_\omega' (x) = 0$ precisely if
$- k(x,l_{-}(x)) l_{-}'(x)
 +
\int_{l_{-}(x)}^{f_\omega (x)} \frac{\partial}{\partial x}  k(x,y) dy = 0$.
\qed

\section{Regularity of solutions of integral equations}\label{sec_ift}


In a number of places in this paper eigenvalue equations $L_a \phi =
\lambda \phi$ for the transfer operator $L_a$ arise. As $L_a$
depends only $C^1$ on $a$, a direct application of the
implicit function theorem as found e.g. in \cite{ber77,chohal82}
yields only weak regularity properties of the solutions. 

The following remark is a variant of Proposition~\ref{prop_smoothnessLa}.
It allows an application of  \cite[Proposition~3.6.1]{buftol03} 
to show that eigenvectors and eigenvalues of $L_a$, in the case of simple eigenvalues,
vary smoothly with $a$.

\begin{remark}\label{rem_toapplybuffoni}
For $a \in I$, $a \mapsto L_a$ 
is a $C^{r+1}$ map from $I$ into ${\cal L}(C^{k+r}(\mathcal{M}),C^k(\mathcal{M}))$, the space of
bounded linear maps from $C^{k+r}(\mathcal{M})$ into $C^k(\mathcal{M})$. 
\end{remark}

We include in this appendix an alternative route 
to obtain such smoothness, as introduction
to the more involved reasoning in Section~\ref{sec_escapeII}.

Given is $L_{a_0} \phi_{a_0} = \lambda_{a_0} \phi_{a_0}$ with
$L_{a}$ acting for $a$ near $a_0$ on $C^k(\mathcal{M})$ 
(here we are considering complex valued functions). 
Similarly we can consider $L_a$ acting on $C^k_0(W)$ for
an isolating neighborhood $W$.
Assume that $\lambda_{a_0}$ is an
isolated eigenvalue of $L_{a_0}$. Denote by $E$ the span of
$\phi_{a_0}$ and let $F^k(W)$ be a complement of $E$ in $C^k(\mathcal{M})$. 
Consider functions $\phi_a = \phi_{a_0} +
\psi_a$ with $\psi_a \in F^k(\mathcal{M})$.
We wish to solve $L_a \phi_a = \lambda_a \phi_a$. Let $P$
be the projection to $E$ along $F^k(\mathcal{M})$. Considering a second
parameter $\lambda$, $L_a \phi_a = \lambda \phi_a$ decomposes as
$$
\left\{ \begin{array}{lcl}
(I-P) L_a (\phi_{a_0}+\psi_a) & = & \lambda  \psi_a,
\\
P L_a  (\phi_{a_0}+\psi_a)& = & \lambda \phi_{a_0}.
\end{array} \right.
$$
The top equation can be solved for $\psi_a$ as a function of $\lambda$ and $a$
for $a$ near $a_0$ and $\lambda$ near $\lambda_{a_0}$.
In fact, by the Fredholm alternative, 
$\psi_a = \left( (I-P)L_a - \lambda I \right)^{-1} \left( (I-P) L_a \phi_{a_0} \right)$.
Putting this into the bottom equation yields a single equation for
$\lambda$. Note that for stationary measures, $\lambda=1$ automatically
solves this equation, compare the proof of Theorem~\ref{continuation}
in Section~\ref{sec_parameter}.

Write the top equation as a fixed point equation
$T_\alpha \psi_\alpha = \psi_\alpha$ with parameters $\alpha$.
The map $T_\alpha$ is a compact linear map mapping $F^k(\mathcal{M})$ into $F^{k+1}(\mathcal{M})$,
compare the proof of Proposition~\ref{compact}. 

\begin{lemma}\label{lemma_plus1}
If $(x,a) \mapsto \psi_a(x)$ is $C^k$, then $(x,a) \mapsto T_a \psi_a (x)$ is $C^{k+1}$.
\end{lemma}

\noindent {\sc Proof.}
See Section~\ref{sec_parameter}.
\qed

\begin{proposition}\label{prop_regular}
Consider the integral equation
$T_\alpha \psi_\alpha = \psi_\alpha$
with $T_\alpha$ as above.
The fixed point $x \mapsto \psi_\alpha(x)$
is smooth jointly in $x,\alpha$.
\end{proposition}

\noindent {\sc Proof.}
Given is a unique fixed point $\psi_\alpha$ depending continuously on $\alpha$.
Formally differentiating $T_\alpha \psi_\alpha = \psi_\alpha$ with respect to $\alpha$ gives
$$
\frac{\partial}{\partial \alpha} \left( T_\alpha \psi_\alpha(x) \right)
=
\frac{\partial}{\partial \alpha} T_\alpha \psi_\alpha(x) + T_\alpha \frac{\partial}{\partial \alpha}
\psi_\alpha (x) = \frac{\partial}{\partial \alpha}\psi_\alpha(x).
$$
 So $\frac{\partial}{\partial \alpha} \psi_\alpha$ should be the solution
$M_\alpha$ of $\frac{\partial}{\partial \alpha} T_\alpha \psi_\alpha (x) + T_\alpha M_\alpha (x)
 = M_\alpha(x)$.
That is,
\begin{equation}\label{BM_a}
M_\alpha  =  (I - T_\alpha)^{-1} \frac{\partial}{\partial \alpha} T_\alpha \psi_\alpha.
\end{equation}
By the Fredholm alternative \cite{kre78}, $I - T_\alpha$ is invertible on $F^k(\mathcal{M})$.
The right hand side of (\ref{BM_a}) is therefore a continuous function.
To establish that $M_\alpha$ is the derivative $\frac{\partial}{\partial \alpha} \psi_\alpha$,
we must show
$|\psi_{\alpha+h} (x) - \psi_\alpha (x) - M_\alpha (x) h| = o(|h|)$ as $h \to 0$
(compare the proof of the implicit function theorem in e.g.
\cite{ber77} or \cite{chohal82}).
Write $\gamma_\alpha (x) = \psi_{\alpha+h} (x) - \psi_\alpha (x)$.
Now
\begin{eqnarray} \nonumber
\gamma_\alpha (x) & = & T_{\alpha+h} (\psi_\alpha + \gamma_\alpha )  (x) - T_\alpha (\psi_\alpha)(x)
\\
\label{Bgammaa}
& = & T_\alpha \gamma_\alpha (x) + \frac{\partial}{\partial \alpha} T_\alpha \psi_\alpha(x) h
+ R(x),
\end{eqnarray}
where
$R(x) = T_{\alpha+h} (\psi_\alpha + \gamma_\alpha )  (x) - T_\alpha (\psi_\alpha)(x)
-  T_\alpha \gamma_\alpha (x) - \frac{\partial}{\partial \alpha} T_\alpha \psi_\alpha(x) h$.
Since $(\psi,\alpha)\mapsto T_\alpha (\psi)$ is differentiable,
for any $\epsilon > 0$ there is $\delta >0$ with
$|R| < \epsilon (|\gamma_\alpha| + |h|)$ if $|\gamma_\alpha|, |h| < \delta$.
Since $\gamma_\alpha$ is continuous in $h$, we may further restrict $\delta$ so that
this estimate on $|R|$ holds for $|h| < \delta$. From
(\ref{Bgammaa}) we get
$\gamma_\alpha = (I - T_\alpha)^{-1} (\frac{\partial}{\partial \alpha} T_\alpha \psi_\alpha h + R)$ so that
$|\gamma_\alpha| < C |h|$ for some $C$, if $|h| < \delta$.
This implies $|R| < \epsilon ( 1 + C) |h|$ for $|h| < \delta$.
As
$$
(I - T_\alpha) (\gamma_\alpha - M_\alpha h)  = R
$$
(from (\ref{BM_a}) we get
$(I - T_\alpha) M_\alpha = \frac{\partial}{\partial \alpha} T_\alpha \psi_\alpha$),
we derive $|\gamma_\alpha (x) - M_\alpha(x) h | < K \epsilon |h|$ for some $K>0$, if $|h| <\delta$.
This proves that $M_\alpha$ equals the
partial derivative $\frac{\partial}{\partial \alpha} \psi_\alpha$.

Higher order derivatives are treated by induction.
Assume that $(x,\alpha) \mapsto \psi_\alpha(x)$
has been shown to be $C^j$.
By Lemma~\ref{lemma_plus1},
$(x,\alpha) \mapsto T_\alpha \psi_\alpha (x)$ is $C^{j+1}$.
So $D \psi_\alpha = D (T_\alpha \psi_\alpha)$ is $C^j$.
As $(I - T_\alpha)^{-1}$ maps $F^j(\mathcal{M})$ to $F^j(\mathcal{M})$, the
right hand side of (\ref{BM_a}) is a $C^j$ function.
The above reasoning shows that $M_\alpha = \frac{\partial}{\partial \alpha} \psi_\alpha$.
Therefore $\frac{\partial}{\partial \alpha} \psi_\alpha$
is $C^j$, so that $(x,\alpha) \mapsto \psi_\alpha(x)$ is $C^{j+1}$.
\qed


\begin{thebibliography}{99}
\footnotesize

\bibitem{abrmarrat83}
R. Abraham, J.E. Marsden, T. Ratiu,
{\em Manifolds, tensor analysis, and applications},
Addison-Wesley Publishing Company, 1983.

\bibitem{alvara03}
J.F. Alves, V. Ara\'ujo,
Random perturbations of nonuniformly expanding maps,
{\em Ast\'erisque} {\bf  286} (2003),  25--62.

\bibitem{alvaravas04}
J.F. Alves, V. Ara\'ujo, C.H. V\'asquez,
 Random perturbations of diffeomorphisms with dominated splitting,
preprint CMUP,  2004.

\bibitem{ara00} V. Ara\'ujo, Attractors and time averages for random maps,
  {\em Ann. Inst. H. Poincar\'e Anal. Non Lin\'eaire} {\bf 17} (2000), 307--369.

\bibitem{arn98}
L.~Arnold, {\em Random Dynamical Systems},
Springer-Verlag, 1998.


\bibitem{bal00}
V. Baladi, {\em Positive transfer operators and decay of correlations},
World Scientific, 2000.


\bibitem{balvia96}
V. Baladi, M. Viana,
Strong stochastic stability and rate of mixing for unimodal maps,
{\em Annales Scientifiques de l'\'Ecole Normale Sup\'erieure S\'er. 4}
{\bf 29}  (1996),  483--517.

\bibitem{balyou93}
V. Baladi, L.-S. Young,
On the spectra of randomly perturbed expanding maps,
{\em Comm. Math. Phys.}  {\bf 156}  (1993),  355--385.

\bibitem{balyou94}
V. Baladi, L.-S. Young,
Erratum: "On the spectra of randomly perturbed expanding maps",
{\em Comm. Math. Phys.} {\bf 166}  (1994),
 219--220.

\bibitem{ber77}
M.S. Berger,
{\em Nonlinearity and functional analysis},
Academic Press, 1977.

\bibitem{bondiavia05} C. Bonatti, L.J. Diaz, M. Viana,
{\em Dynamics beyond uniform hyperbolicity}, Springer-Verlag, 2005.

\bibitem{braolivar98}
S. Brassesco, E. Olivieri, M.E. Vares,
Couplings and asymptotic exponentiality of exit times,
{\em J. Statist. Phys.} {\bf 93} (1998),  393--404.

\bibitem{buftol03}
B. Buffoni, J. Toland, {\em Analytic theory of global
bifurcation}, Princeton University Press, 2003.

\bibitem{chemartro98}
N. Chernov, R.  Markarian, S. Troubetzkoy,
Conditionally invariant measures for Anosov maps with small holes,
{\em Ergodic Theory Dynam. Systems} {\bf  18}  (1998),  1049--1073.

\bibitem{chohal82}
S.-N. Chow, J.K. Hale,
{\em Methods of bifurcation theory}, Springer-Verlag, 1982.

\bibitem{colkli00}
F. Colonius, W. Kliemann,
{\em The dynamics of control},
Birkh\"auser,  2000.

\bibitem{colgaykli04}
F. Colonius, T. Gayer, W. Kliemann,
Near invariance for Markov diffusion systems,  preprint.

\bibitem{day83}
M.V. Day,
On the exponential exit law in the small parameter exit problem,
{\em Stochastics} {\bf 8} (1983),  297--323.



\bibitem{delljun99} M. Dellnitz, O. Junge, On the approximation of complicated
dynamical behavior, {\em SIAM J. Numer. Anal.} {\bf 36} (1999), 491--515.

\bibitem{doiinosad98}
S. Doi, J. Inoue, S. Kumagai,
Spectral analysis of stochastic phase lockings and
stochastic bifurcations in the sinusoidally forced
van der Pol oscillator with additive noise,
{\em  J. Statist. Phys.}  {\bf 90}  (1998),  1107--1127.

\bibitem{doo53} J.L. Doob, {\em Stochastic Processes}, John Wiley \& Sons, 1953.

\bibitem{eckthowit81}
 J.-P. Eckmann, L.  Thomas, P. Wittwer,
Intermittency in the presence of noise,  {\em J. Phys. A}  {\bf 14}  (1981), 3153--3168.


\bibitem{ferkesmarpic95}
P.A. Ferrari, H. Kesten, S. Martinez, P. Picco,
Existence of quasi-stationary distributions. A renewal dynamical approach,
{\em Ann. Probab.} {\bf 23} (1995), 501--521.

\bibitem{for51}
M.K. Fort,
Points of continuity of semi-continuous functions,
{\em Publ. Math. Debrecen} {\bf  2} (1951), 100--102.

\bibitem{frewen84}
M.I. Freidlin, A.D. Wentzell,
{\em Random perturbations of dynamical systems},
Springer-Verlag,  1984.

\bibitem{gay01}
T. Gayer,
On Markov chains and the spectra of the corresponding Frobenius-Perron operators,
{\em Stoch. Dyn.} {\bf 1} (2001),  477--491.

\bibitem{gay04}
T. Gayer,
Control sets and their boundaries under parameter variation,
{\em J. Differential Equations} {\bf 201} (2004), 177--200.


\bibitem{greottromyor87}
C. Grebogi, E. Ott, F.  Romeiras, J.A.  Yorke,
Critical exponents for crisis-induced intermittency, {\em  Phys. Rev. A}  {\bf 36}  (1987),
5365--5380.

\bibitem{her79}
M.R. Herman, Sur la conjugaison diff\'erentiable des diff\'eomorphismes
du cercle \`a des rotations,
{\em Publications math\'ematique de l'I.H.\'E.S.} {\bf 49} (1979), 5--233.


\bibitem{hirhubsca82}
J.E. Hirsch, B.A. Huberman, D.J. Scalapino,
Theory of intermittency,
{\em Physical Review A} {\bf  25} (1982),  519--532.

\bibitem{hir76}
 M.W. Hirsch, {\em Differential topology}, Springer-Verlag, 1976.

\bibitem{homyou02}
A.J. Homburg, T. Young,
Intermittency in families of unimodal maps,
{\em Ergodic Theory Dynam. Systems} {\bf 22} (2002),  203--225.

\bibitem{homyou05}
A.J. Homburg, T. Young,
Intermittency and Jakobson's theorem near saddle-node bifurcations, preprint.

\bibitem{homyou05b}
A.J. Homburg, T. Young,
Hard bifurcations in dynamical systems with bounded random perturbations,
preprint.

\bibitem{kat66} T. Kato, {\em Perturbation theory for linear operators},
Springer-Verlag, 1966.

\bibitem{kif86} Y. Kifer, {\em Ergodic theory of random transformations}, Birkh\"auser, 1986.

\bibitem{kif88} Y. Kifer, {\em Random perturbations of dynamical systems}, Birkh\"auser, 1988.

\bibitem{kre78}
E. Kreyszig, {\em Introductory functional analysis with applications},
John Wiley \& Sons, 1978.



\bibitem{lasmac94}
A. Lasota, M.C. Mackey, {\em Chaos, fractals and noise}, Springer-Verlag, 1994.

\bibitem{laspea01}
J.B. Lasserre, C.E.M. Pearce,
On the existence of a quasistationary measure for a Markov chain,
{\em Ann. Probab.} {\bf 29}, 427--446.


\bibitem{mac91} R.S. MacKay, An extension of Zeeman's notion of structural stability
to noninvertible maps, {\em Phys. D} {\bf 52} (1991), 246--253.

\bibitem{melstr93}
W. de Melo, S. van Strien,
{\em One-dimensional dynamics}, Springer-Verlag, 1993.

\bibitem{mil00}
J.W. Milnor, {\em Dynamics: introductory lectures}, \verb#http://www.math.sunysb.edu/~jack/DYNOTES/#

\bibitem{pia81}
G. Pianigiani,
 Conditionally invariant measures and exponential decay,
{\em J. Math. Anal. Appl.}  {\bf  82}  (1981),  75--88.

\bibitem{piayor79}
G. Pianigiani, J.A. Yorke,
Expanding maps on sets which are almost invariant. Decay and chaos,
{\em  Trans. Amer. Math. Soc.} {\bf  252}  (1979), 351--366.

\bibitem{pomman80}
Y. Pomeau, P. Manneville,
Intermittent transition to turbulence in dissipative dynamical systems,
{\em Comm. Math. Phys.} {\bf 74} (1980),  189--197.

\bibitem{qua91} A.N. Quas, On representations of Markov chains by random smooth maps,
{\em Bull. London Math. Soc.} {\bf 23} (1991), 487--492.

\bibitem{rue82}
D. Ruelle,
Small random perturbations of dynamical systems
and the definition of attractors,  {\em Comm. Math. Phys.}  {\bf 82}  (1982),  137--151.

\bibitem{rue90}
D. Ruelle, An extension of the theory of Fredholm determinants,
{\em Inst. Hautes \'Etudes Sci. Publ.Math.} {\bf 72} (1990), 175--193.

\bibitem{sch74}
H.H. Schaefer, {\em Banach lattices and positive operators},
Springer-Verlag, 1974.


\bibitem{via97} M. Viana, {\em Stochastic dynamics of deterministic systems},
Col. Bras. de Matem\'atica, 1997.

\bibitem{yos80} K. Yosida, {\em Functional Analysis}, Springer-Verlag, 1980.

\bibitem{you86}
L.-S. Young,
Stochastic stability of hyperbolic attractors,
{\em Ergodic Theory Dynam. Systems} {\bf  6} (1986),  311--319.

\bibitem{zee88} E.C. Zeeman, Stability of dynamical systems,
{\em Nonlinearity} {\bf 1} (1988),
115--155.


\end{thebibliography}
\end{document}